\documentclass[12pt]{amsart}

\usepackage{amsmath}

\usepackage{amssymb}
\usepackage{amsthm}
\usepackage{psfig}
\newcommand{\ignore}[1]{\relax}

\newcommand{\C}{\mathbb C}
\newcommand{\R}{\mathbb R}
\newcommand{\Z}{\mathbb Z}

\newcommand{\NN}{\mathcal N}
\newcommand{\Q}{\mathbb Q}
\newcommand{\U}{\mathcal U}
\newcommand{\F}{\mathcal F}
\newcommand{\DD}{\mathcal D}
\newcommand{\HH}{\mathcal H}
\newcommand{\LL}{\mathcal Log}

\newcommand{\V}{{V}^{\circ}}
\newcommand{\ve}{\stackrel{\to}{v}}
\newcommand{\PP}{\mathcal P}

\newcommand{\const}{\operatorname{const}}
\newcommand{\val}{\operatorname{val}}

\newtheorem{thm}{Theorem}
\newtheorem{thmquote}{Theorem}

\newtheorem{lem}{Lemma}[section]
\newtheorem{cor}[lem]{Corollary}
\newtheorem{prop}[lem]{Proposition}

\newtheorem{que}{Question}

\theoremstyle{definition}
\newtheorem{defn}{Definition}
\newtheorem{exa}{Example}
\theoremstyle{remark}
\newtheorem{rem}[lem]{Remark}

\newtheorem{rmk}[lem]{Remark}
\newcommand{\tor}{(\C^*)^{n+1}}
\newcommand{\ktor}{(K^*)^{n+1}}
\newcommand{\rtor}{(\R^*)^{n+1}}

\newcommand{\conj}{\operatorname{conj}}
\newcommand{\dd}{\partial}
\newcommand{\am}{\mathcal{A}}
\newcommand{\D}{\mathcal{D}}
\newcommand{\cp}{{\mathbb C}{\mathbb P}}
\newcommand{\rp}{{\mathbb R}{\mathbb P}}

\newcommand{\Log}{\operatorname{Log}}

\newcommand{\Vol}{\operatorname{Vol}}

\newcommand{\Int}{\operatorname{Int}}
\renewcommand{\setminus}{\smallsetminus}

\begin{document}

\title
[Decomposition into higher-dimensional pairs-of-pants]
{Decomposition into
pairs-of-pants
for complex algebraic hypersurfaces
}
\author{Grigory Mikhalkin}
\address {Department of Mathematics\\
University of Utah\\
Salt Lake City, UT 84112, USA}
\address{St.Petersburg Branch, Steklov Mathematical Institute,
Fontanka 27, St.Petersburg, 191011, Russia}
\email{mikhalkin@math.utah.edu}

\begin{abstract}
It is well-known that a Riemann surface can be decomposed
into the so-called {\em pairs-of-pants}.
Each pair-of-pants is diffeomorphic to
a Riemann sphere minus 3 points.
We show that a smooth complex projective hypersurface of
arbitrary dimension admits a
similar decomposition. The $n$-dimensional pair-of-pants
is diffeomorphic to $\cp^n$ minus $n+2$ hyperplanes.

Alternatively, these decompositions can be treated
as certain fibrations on the hypersurfaces.
We show that there exists a singular fibration on the hypersurface
with an $n$-dimensional polyhedral complex as its base
and a real $n$-torus as its fiber. The base accommodates
the geometric genus of a hypersurface $V$.
Its homotopy type is a wedge of $h^{n,o}(V)$ spheres $S^n$.
\end{abstract}

\maketitle

\section{Introduction}
\subsection{Main question}
In this paper we study non-singular algebraic hypersurfaces
in $\cp^{n+1}$ and other toric varieties.
Let $V$ be such a hypersurface.
Naturally, $V$ is a complex variety and thus
has the underlying structure of a smooth manifold.
Furthermore, $V$ is a symplectic manifold.
The symplectic structure is induced by the embedding to $\cp^{n+1}$.

Since $V$ is non-singular, its diffeomorphism and symplectomorphism
types depend only on its degree,
i.e. the degree of the defining polynomial $f$.
All smooth hypersurfaces of the same degree are isotopic
in the ambient $\cp^{n+1}$ even though the complex structure of $V$ varies
with the coefficients of $f$.

Thus, from the point of view of differential topology or
symplectic topology a smooth projective hypersurface $V$ is
given by two numbers: its dimension $n$ and its degree $d$.

\begin{que}
Given $n$ and $d$, describe a non-singular hypersurface
$V\subset\cp^{n+1}$ of degree $d$
as a smooth manifold and as a symplectic manifold.
\end{que}

More generally, one can ask a similar question where $\cp^{n+1}$
is replaced by an arbitrary toric variety.
The degree $d$ would then
be replaced with a convex lattice polygon $\Delta\subset\R^{n+1}$.

\subsection{State of knowledge for small values of $n$ and $d$}
\subsubsection{Case $n=1$}
The answer to this question is well-known if $n=1$. Then $V$
is a Riemann surface. Topologically it is a sphere with $g$
handles, where the genus $g$ can be computed from the degree $d$
by the adjunction formula $g=\frac{(d-1)(d-2)}{2}$.

Recall that one way to understand Riemann surfaces
is via their decomposition to primitive pieces each diffeomorphic
to a sphere with 3 holes.
These primitive pieces are called {\em pairs-of-pants}
and such a decomposition can be thought of as some (singular)
fibration of the Riemann surface over a 3-valent graph,
see Figure \ref{tripod}.
Note that the first Betti number of the base graph coincides
with the genus of the Riemann surface.

\subsubsection{Case $n=2$}
In this case $V$ is a smooth 4-manifold.
If $d$ is 1, 2 or 3 then $V$ is diffeomorphic to $\cp^2$,
$\cp^1\times \cp^1$ or $\cp^2\#6\bar{\cp}^2$, the connected sum
of $\cp^2$ and 6 copies of $\cp^2$ with the inverse orientation.
In these cases the geometric genus $p_g=h^{2,0}(V)$ vanishes.

If $d=4$ then $p_g=1$ and $V$ is the celebrated K3 surface (named so,
according to A. Weil, in honor of K\"ahler, Kodaira, Kummer and
the K2-mountain in Pakistan). This manifold is primitive,
it does not decompose
as a connected sum.
One way to understand its topology
is via a singular fibration $\lambda:V\to S^2$. A generic
fiber of $\lambda$ is a torus while 24 fibers are special
and are homeomorphic to a torus with its meridian collapsed
to a point (so-called fishtail fibers). The fibration $\lambda$
can be chosen so that all generic fibers are Lagrangian submanifolds,
i.e. so that the symplectic 2-form restricted to these fibers vanishes.

For any value of $d$ $V$ is simply-connected and if $d\ge 4$
it does not decompose into a connected sum. We have
$p_g(V)=\frac{(d-1)(d-2)(d-3)}{6}$ (cf. e.g. a more
general Khovanskii's formula \cite{Kh}).
A simply-connected smooth 4-manifolds is determined
up to a homeomorphism
once we know its Euler characteristic $\chi$, its signature $\sigma$
and whether it is spin or not.
Our manifold $V$ is spin iff $d$ is even, $\chi=d^3-4d^2+6d$
and $\sigma=2(2p_g+1)-(\chi-2)=4(p_g+1)-\chi=
\frac{4d-d^3}{3}$.

The diffeomorphism (and symplectomorphism) type of $V$ is, however,
more mysterious as it is not determined by purely homological data.
E.g. the surface of degree 5 is a non-spin manifold with $\chi=55$
and $\sigma=-35$, but there might be many non-diffeomorphic
manifolds with these data.

\subsubsection{Case $d=n+2$}
In this case the canonical class of $V$ is
trivial and there exists a nowhere-degenerate holomorphic
n-form $\Omega$ on $V$.
Such $V$ is called a {\em Calabi-Yau} manifold.
Here we have $p_g=1$.
According to the Strominger-Yau-Zaslow conjecture \cite{SYZ} there
is supposed to exist a singular special Lagrangian fibration of $V$
over the sphere $S^n$. This means that a generic fiber should
be Lagrangian and such that the imaginary part of $\Omega$
restricted to the fiber is zero as a real 3-form at every point.

It was verified in \cite{Zh} and \cite{Ru} that such fibrations
exist in this case at least if we relax a special Lagrangian condition
to simply Lagrangian.
Note that special Lagrangian condition makes use of
the non-degenerate holomorphic $n$-form from a Calabi-Yau manifold.
Thus, at least literally, the Strominger-Yau-Zaslow conjecture
only makes sense if $d=n+2$ in our setup.
However, a relaxed version of this conjecture makes sense
for all values of $d$ and $n$.

\subsection{Results of the paper}
Here we state the main results of the paper informally.
See section \ref{mainres} for precise statements.

\subsubsection{Torus fibration and pairs-of-pants decomposition}
Theorem \ref{diff} asserts that for any value of $n$ and $d$
the hypersurface $V$ admits a singular fibration $\lambda$ over
an $n$-dimensional polyhedral complex $\bar\Pi$.
A generic fiber of $\lambda$ is diffeomorphic to
a smooth torus $T^n$.
The base $\bar\Pi$ here is homotopy equivalent to the bouquet
of $p_g$ copies of $S^n$, thus this theorem can be interpreted
as a geometric interpretation of the geometric genus $p_g$.

Furthermore, the local topological structure of
the polyhedral complex $\Pi\subset\R^{n+1}$ is
known in differential topology as the local structure
of so-called {\em special spines}. In particular,
there is a natural stratification of $\Pi$
and regular neighborhoods of the vertices essentially
exhaust the complex $\Pi$.

It turns out that the stratification of the base $\Pi$
determines a decomposition of the hypersurface $V$ into
$d^{n+1}$ copies of $\PP_n$, where $\PP_n$ is diffeomorphic
to $\cp^n$ minus $(n+2)$ hyperplanes in general position.
This decomposition can be considered as a higher-dimensional
analogue of the pair-of-pants decomposition of Riemann surfaces.
In particular $\PP_1$ is the classical pair-of-pants $\hat{\C}\setminus
\{0,1,\infty\}$.

\subsubsection{A projective hypersurface as a piecewise-linear object}
The base $\Pi$ of the fibration $\lambda$
is a piecewise-linear $n$-dimensional complex in $\R^{n+1}$.
The dimension over $\R$ of the hypersurface $V$ is $2n$.
Yet the hypersurface $V$ can be reconstructed (as a smooth manifold)
from $\Pi\subset\R^{n+1}$. It turns out that $\Pi$ (together
with its PL-embedding to $\R^{n+1}$) encodes the combinatorics
of gluing of $d^{n+1}$ copies of $\PP_n$ needed to obtain $V$.
Theorem \ref{locthm} is the corresponding reconstruction theorem.

\subsubsection{Lagrangian submanifolds in projective hypersurfaces}
It turns out that the fibration $\lambda$ produces a number
of Lagrangian submanifolds in $V$.
Different fibers of $\lambda$ are not necessarily homologous
and $p_g=h^{n,0}$ disjoint embedded Lagrangian tori
come as fibers of $\lambda$. This tori are linearly independent
in $H_n(V)$. In addition we have $h^{n,0}$ linearly independent
embedded Lagrangian spheres coming as partial sections of $\lambda$.
In particular, we have Corollary \ref{lagr}.

\subsection{Acknowledgement}
A large part of this paper was written during author's visit
to Rio de Janeiro in February 2002.
The author thanks Instituto de Matem\'atica Pura e Aplicada
for providing outstanding conditions for research and writing.
The author is supported in part by the NSF grant DMS-0104727.

\section{Preliminaries}
\subsection{Balanced polyhedra}
\begin{defn}
\label{prpc}
A subset $\Pi\subset\R^{n+1}$ is called a
{\em proper rational polyhedral complex} (or just
a {\em polyhedral complex} in this paper) if it can be
presented as a finite union of closed sets in $\R^{n+1}$
called {\em cells} with the following properties.
\begin{itemize}
\item Each cell is a closed convex (possibly semi-infinite)
polyhedron. The dimension of the cell is, by definition,
the dimension of its affine span,
the smallest affine subspace of $\R^{n+1}$ which contains it.
We call a cell of dimension $k$ a $k$-cell.
\item The slope of the affine span of each cell is rational.
I.e. the linear subspace of $\R^{n+1}$ parallel to the affine
span is defined over $\Q$.
\item The boundary (i.e. the boundary in the corresponding affine span)
of a $k$-cell is a union of $(k-1)$-cells.
\item Different open cells (i.e. the interiors of the cells in
the corresponding affine spans) do not intersect.
\end{itemize}
\end{defn}

Informally speaking, a proper polyhedral complex in $\R^{n+1}$
is a cellular space where each cell is a convex polyhedron
with a rational slope and
where some cells are allowed to go to infinity.

As usual, the dimension of $\Pi$ is the maximal
dimension of its cells.

\begin{defn}
A polyhedral $n$-complex is called {\em weighted} if there is
a natural number $w(F)$, called {\em weight}, prescribed to each of
its $n$-cell $F$. (Of course, any polyhedral complex can be considered
as a weighted polyhedral complex by prescribing 1 to each $n$-cell.)
\end{defn}

Let $\Pi\subset\R^{n+1}$ be a weighted polyhedral $n$-complex.
Note that its complement $\R^{n+1}\setminus\Pi$ consists of
a finite union of connected components.
Let $F\subset\Pi$ be an $n$-cell.

Recall that by Definition \ref{prpc} the $n$-cell $F$
has a rational slope in $\R^{n+1}$. Therefore, it defines
an integer covector $$\pm c_F:\Z^{n+1}\to\Z$$ up to its sign.
Here are the characteristic properties of $c_F$.
\begin{itemize}
\item The kernel of $c_F$ is parallel to $F$.
\item $\frac{1}{w(F)}c_F$ is a primitive (i.e. non-divisible)
integer covector
$\Z^{n+1}\to\Z$.
\end{itemize}

Furthermore, even the sign of $c_F$ becomes well-defined
once we co-orient $F\subset\R^{n+1}$.

Polyhedral complexes
that appear in this paper have the following additional property.

\begin{defn}
A weighted polyhedral $n$-complex $\Pi\subset\R^{n+1}$
is called {\em balanced} if for every $(n-1)$-cell $G\subset\Pi$
the following condition holds.
Let $F_1,\dots,F_k$ be the $n$-cells adjacent to $G$.
A choice of a rotational direction about $G$ defines
a coherent co-orientation on these $n$-cells.
The balancing condition is
$$\sum\limits_{j=1}^k c_{F_j}=0.$$
\end{defn}

\begin{figure}[h]
\centerline{\psfig{figure=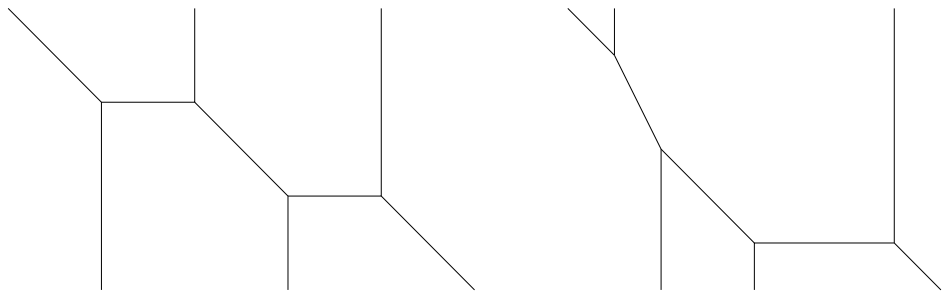,height=1.5in,width=4.9in}}
\caption{\label{ball} Balanced graphs in $\R^2$.}
\end{figure}

\begin{exa}
\label{si}
Consider the function
$$H(x_1,\dots,x_{n+1})=\max\{0,x_1,\dots,x_{n+1}\}.$$
This is a convex piecewise-linear function $\R^{n+1}\to\R$.
We define the {\em primitive complex} $\Sigma_n\subset\R^{n+1}$
as the {\em corner locus} of $H$, i.e. the set of points where
$H$ is not smooth.

Note that $\Sigma_n$ is a balanced proper polyhedral
complex in $\R^{n+1}$. Its $k$-cells are formed by
the points where at least $n+2-k$ of the functions
$0,x_1,\dots,x_{n+1}$ achieve the value of $H$.
In fact, it is easy to see that topologically $\Sigma_n$
is the cone over the $(n-1)$-skeleton of the $(n+1)$-simplex.
The fact that $\Sigma$ is balanced follows from Proposition \ref{bal}.
\end{exa}
\begin{figure}[h]
\centerline{\psfig{figure=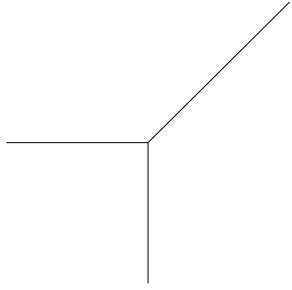,height=1.5in,width=1.5in}}
\caption{Primitive complex $\Sigma_n$.}
\end{figure}

The following example is a generalization of the previous one.
As the following propositions show, it is the fundamental example
of balanced polyhedra.

\begin{exa}
\label{tri}
Let $A\subset\Z^{n+1}$ be a finite set and
let $v:A\to\R$ be any function. Let $\Delta\subset\R^{n+1}$
be the convex hull of $A$.
We associate the following polyhedral complex $\Pi_v$ to $v$.

Take the Legendre transform of $v$, $L_v:\R^{n+1}\to\R,$
$$L_v(y)=\max\limits_{x\in A} (xy-v(x)).$$
Here $x,y\in\R^{n+1}$ and $xy$ is their scalar product.
Since the maximum is taken over a finite set,
the result $L_v$ is a convex piecewise-linear function.
We define $\Pi_v$ as the corner locus of $L_v$
(recall that this is the set of points where $L_v$ is not smooth).

To present Example \ref{si} as a special case of Example \ref{tri}
we take the vertices of the standard simplex
\begin{equation}
\label{sta}
\Delta_1 \{(x_1,\dots,x_{n+1})\in\R^{n+1}\ |\ x_j\ge 0,
x_1+\dots+x_{n+1}\le 1 \}
\end{equation}
for $A$ and set $v\equiv 0$.
\end{exa}

Recall that a polyhedron in $\R^{n+1}$ is called {\em a lattice
polyhedron} if all its vertices belong to $\Z^{n+1}$. A
subdivision of a polyhedron into smaller polyhedra is called {\em
a lattice subdivision} if all its subpolyhedra are lattice.

\begin{prop}
\label{podrazd}
The set $\Pi_v$ from Example \ref{tri} is a proper
rational polyhedral complex dual to a certain lattice
subdivision of $\Delta$.
\end{prop}
\begin{proof}
We start by associating to $v$ a certain lattice subdivision
${\mathcal D}_v$ of $\Delta$.
Let $O\Gamma(v)$ be the overgraph of $v$,
i.e. the set of vertical rays upwards in $\R^{n+1}\times\R$
starting at the points of the graph of $v$.
The convex hull of $O\Gamma(v)$ is a semi-infinite
closed polyhedral domain. The projections of its finite
faces to $\R^{n+1}$ form the subdivision ${\mathcal D}_v$.

We claim that $\Pi_v$ is a polyhedral
complex dual to ${\mathcal D}_v$.
Namely, a $k$-dimensional polyhedron $\Delta'$ in ${\mathcal D}_v$, $k>0$,
gives a $(n+1-k)$-cell of $\Pi_v$.
This cell is compact iff $\Delta'\subset\Delta$.

This claim follows from the duality property of the Legendre
transform. Consider the function $\underline{v}$ whose graph
is is given by the lower boundary of the convex hull of $O\Gamma(v)$.
If $v$ is convex then the function $\underline{v}$ extends $v$ and
is defined on the whole
polyhedron $\Delta$, not just on its lattice points.
It is a convex piecewise-linear function. The Legendre transform
of $v$ coincides with the Legendre transform of $\underline{v}$.
(In fact the function $\underline{v}$ can be defined by applying
the Legendre transform to $v$ twice.)
By duality, the graph of $L_{\underline{v}}$ has the facets en lieu
of the vertices of the graph of $\underline{v}$ and so on.
\end{proof}

Note that $\Pi_v$ is naturally weighted.
Indeed, an $n$-cell $F\subset\Pi_v$ comes as a corner
between the graphs of two integer linear functions.
The difference between these functions is an integer
covector $c_F$. We define $w(F)\in{\mathbb N}$
as the maximum integer divisor of $c_F$.

\begin{prop}
\label{bal}
The weighted polyhedral complex $\Pi_v$ is balanced.
\end{prop}
\begin{proof}
The proposition easily follows from the definition of the
covectors $c_{F_j}$ for the $n$-cells $F_j$ adjacent to
an $(n-1)$-cell $G\subset\Pi$.
\end{proof}


%
%

\begin{rem}
\label{ambig}
Note that several different functions $v$ define the same
complex $\Pi_v$ by the construction of Example \ref{tri}.
Here is the list of ambiguities.
\begin{enumerate}
\item Let $v'=v+\const:A\to\R$ be a function different
with $v$ by a constant. Then $\Pi_v=\Pi_{v'}$.
\item Let $A'=A+c$, where $c\in\Z^{n+1}$ and
$v':A'\to\R$ is defined by $v'(z+c)=v(z)$.
Then $\Pi_v=\Pi_{v'}$.
\item Let $A'$ be such that its convex hull $\Delta'$
coincides with $\Delta$, the convex hull of $A$.
Let $\underbar{v}$ (resp. $\underbar{v'}$)
be the maximal convex function such that $\underbar{v}\le v$
(resp. $\underbar{v'}\le v'$). Suppose that
$\underbar{v}=\underbar{v'}$. Then $\Pi_v=\Pi_{v'}$.
\end{enumerate}
\end{rem}

The following proposition shows that Example \ref{tri}
is fundamental.
\begin{prop}
\label{obratno}
Suppose that $\Pi\subset\R^{n+1}$ is a weighted balanced
proper rational polyhedral complex. Then there exists
a finite set $A\subset\Z^{n+1}$ and a function $v:A\to\Z$ such
that $\Pi=\Pi_v$ (see Example \ref{tri}). The convex
hull $\Delta\subset\R^{n+1}$ of $A$ is unique up to a translation
in $\Z^{n+1}$. The choice of the
function $v$ is unique up to the ambiguity of Remark \ref{ambig}.
\end{prop}
\begin{proof}
First we define a convex piecewise-linear function $H$
whose corner locus is $\Pi$ and then choose a function
$v$ such that $H$ is the Legendre transform $L_v$ of $v$.
Note that the finiteness condition in Definition \ref{prpc}
implies that there are finitely many connected components
in $\R^{n+1}\setminus\Pi$.

We define the function $H$ inductively.
Choose any connected component $D_0$ of $\R^{n+1}\setminus\Pi$ as
a ``reference component". Define $H|_{D_0}\equiv 0$.
Suppose that $D'$ is a component of $\R^{n+1}\setminus\Pi$
such that there exists an adjacent component $D$ where $H$
is already defined.

Let $F$ be the $n$-cell of of $\Pi$
separating $D$ from $D'$. Let $c_F$ be the covector associated
to $F$ (recall that the weight of $F$ is incorporated into $c_F$)
with the co-orientation directed from $D$ to $D'$.
Let $l_D:\R^{n+1}\to\R$ be the affine-linear function extending
$H|_D$. We define $H|_{D'}=l_D+c_F+c$, where the constant $c$
is chosen so that $H|_D$ and $H|_{D'}$ agree on $F$.
By the balancing
condition the result does not depend on the choice of
the adjacent component $D$ where $H$ is already defined.

To define $v$ we take the Legendre transform of $H$.
This amounts to associating to each component $D$ a point
$z\in\Z^{n+1}$ equal to the gradient of $H|_D$ and
setting $v(z)=-l_D(0)$. Thus, the number of elements
of the set $A$ is equal to the number of components of
$\R^{n+1}\setminus\Pi$.

The ambiguity Remark \ref{ambig}.3 comes from taking
the Legendre transform of non-convex functions $v$.
It coincides with the Legendre transform of the
underlying convex function $\underbar{v}$.
(In fact, nothing changes if we assume that $v$
is defined on the whole $\Z^{n+1}$ by letting
$v(z)=+\infty$ for $z\notin A$.)
The ambiguities Remark \ref{ambig}.1 and \ref{ambig}.2
come from the ambiguity in assigning a linear function
for $H|_{D_0}$.
\end{proof}

\begin{cor}
\label{dual}
Any $n$-dimensional balanced
polyhedral complex $\Pi\subset\R^{n+1}$ determines
a convex lattice polyhedron $\Delta\subset\R^{n+1}$
(defined up to translation) and a lattice subdivision
of $\Delta$.
\end{cor}
This corollary follows from Propositions \ref{obratno} and
\ref{podrazd}.
\begin{figure}[h]
\centerline{\psfig{figure=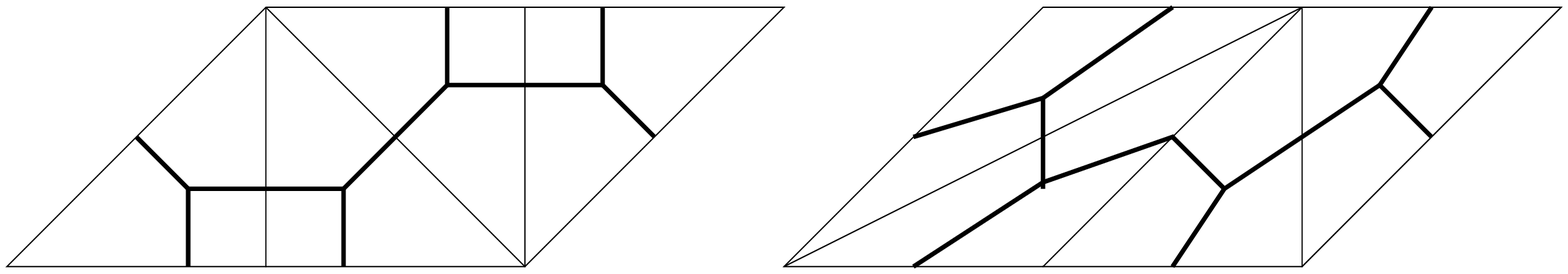,height=0.85in,width=4.9in}}
\caption{\label{podr} The lattice polyhedron subdivisions dual to the
balanced graphs from Figure \ref{ball}.}
\end{figure}

The next corollary illustrates the strength
of the balancing condition that we require just at the $n$-cells.
We do not use this corollary elsewhere in the paper.

Let $B$ be a vertex of $\Pi$ and let $E_1,\dots,E_k$ be
the edges adjacent to $B$. Let $v_j\in\Z^{n+1}$, $j=1,\dots,k$,
be the primitive integer vectors parallel to $E_j$
and directed outwards from $B$.
Suppose that each $E_j$ is adjacent to exactly $n+2$
connected components of $\R^{n+1}\setminus\Pi$
(note that this is a general position situation).
\begin{cor}
If $\Pi\subset\R^{n+1}$ is a balanced $n$-complex
then there exists a weight $w_j\subset{\mathbb N}$ for $E_j$,
$j=1,\dots k$, such that
$$\sum\limits_{j=1}^k w_jv_j=0.$$
\end{cor}
\begin{proof}
By Proposition \ref{obratno} $\Pi$ comes as a corner
locus of a convex piecewise-linear function $F$ on $\R^{n+1}$.
Let $y=a_{j,1}x_1+\dots+a_{j,n+1}x_{n+1}$, $j=1,\dots,n+2$
be the equations of the linear functions on the adjacent
components of $\R^{n+1}\setminus\Pi$.
Then $u_j=(a_{j,1},\dots,a_{j,n+1},-1)$ are
the vectors in $\R^{n+2}=\R^{n+1}\times\R$ normal
to the linear portions of the graph of $F$ adjacent to $B$.

The $\R^{n+2}$-version of the vector product associates
a normal vector to $(n+1)$ other vectors in $\R^{n+2}=\R^{n+1}\times\R$.
We take all possible such products among $u_j$ and project them
to $\R^{n+1}$. The result is the vectors which are multiples of $v_j$.
By linear algebra the sum of these vectors is zero.
\end{proof}

\subsection{Maximal polyhedral complexes and their decomposition
into primitive pieces}
\label{s1.2}
\begin{defn}
We call $\Pi$ a {\em dual $\Delta$-complex} if it corresponds
to the convex polyhedron $\Delta\subset\R^{n+1}$
by Proposition \ref{obratno}.
We call $\Pi$ a {\em maximal} polyhedral complex
if the elements of
the corresponding subdivision from Corollary \ref{dual}
are simplices of volume $\frac{1}{(n+1)!}$ (a so-called
{\em unimodular} lattice triangulation).
\end{defn}

\begin{prop}
The minimal positive volume of a lattice polyhedron
in $\R^{n+1}$ is $\frac{1}{(n+1)!}$.
Any lattice polyhedron of volume $\frac{1}{(n+1)!}$
can be identified with the standard simplex $\Delta_1$
(see \eqref{sta}) by an element of
$ASL_{n+1}(\Z)$.
\end{prop}
Here $ASL_{n+1}(\Z)$
stands for the group of affine-linear transformations
of $\R^{n+1}$ whose rotation part belongs
to $SL_{n+1}(\Z)$.
\begin{proof}
We may assume that our lattice polyhedron is a simplex,
since otherwise we can triangulate it to smaller polyhedra.
Fix one of its vertice  and consider the $(n+1)$ integer
vectors connecting it to other vertices.
The volume of the simplex is equal to the determinant of the
sublattice generated by these vectors divided by $(n+1)!$.
\end{proof}

\begin{exa}
Clearly, a dual $\Delta_1$-complex (see \eqref{sta})
is necessarily maximal.
The complexes from Figure \ref{ball} are
maximal dual $\Delta$-complexes for the polyhedra $\Delta$
pictured on Figure \ref{podr}.
\end{exa}

\begin{prop}
\label{prim}
Any dual $\Delta_1$-complex is the result of a translation of
$\Sigma_n$ in $\R^{n+1}$.
\end{prop}
\begin{proof}
Such a complex $\Pi$ is determined by a function
$v:\Delta_1\cap\Z^{n+1}\to\R$, i.e. by $n+2$ numbers
$a_1,\dots,a_{n+1},b\in\R$.
Recall (see Example \ref{tri}) that $\Pi$ is the corner
locus $L_v(x_1,\dots,x_{n+1})=\max\{x_j-a_j,-b\}$.
If $a_j=b=0$ for all $j$ than $\Pi=\Sigma_1$.
Adding the same real number to all numbers does not change $\Pi$.
Changing $a_j$ by $t$ results in a translation by $t$ in the direction
of $x_j$.
\end{proof}


%
\begin{rmk}
\label{nemax}
Not for every $\Delta$ there exists maximal dual $\Delta$-complex.
E.g. a lattice simplex in $\R^3$, whose vertices
are $(1,0,0)$, $(0,1,0)$, $(1,1,0)$ and $(0,0,n)$,
cannot be further subdivided.
On the other hand, a maximal dual $\Delta$-complex is, of course,
not unique.
\end{rmk}

\begin{prop}
\label{genus}
If $\Pi$ is a maximal dual $\Delta$-complex then
$\Pi$ is
homotopy equivalent to the bouquet of
$\#(\Int\Delta\cap\Z^{n+1})$ copies of $S^n$.
\end{prop}
\begin{proof}
Because of its maximality, the polyhedron
$\Pi$ is dual to a unimodular triangulation
of $\Pi$. Such a triangulation cannot be further subdivided and
therefore its vertices are all the lattice points of $\Delta$.
Therefore, $\Pi$ is homotopy equivalent to
$\Int\Delta\setminus\Z^{n+1}$.
\end{proof}

Here is a way to canonically cut a maximal complex $\Pi\subset\R^{n+1}$
into standard-looking subsets $U_j$.
We define the cutting locus $\Xi$ as the following simplicial
complex that is partially dual to $\Pi$.
The vertices of $\Xi$ are the baricenters of all bounded
$k$-cells, $k>0$ from $\Pi$.
The simplices of $\Xi$ have the baricenters of positive-dimensional
cells $F_k\subset\Pi$ in the embedded
towers $F_1\subset\dots\subset F_l$ as its vertices.
Note that $\Xi\subset\Pi$ is a finite simplicial $(n-1)$-complex.

%

\begin{defn}
\label{uj}
The connected components of $\Pi\setminus\Xi$
are called {\em the primitive pieces} of $\Pi$.
We denote them with $U_j$.
These open sets are parametrized by the vertices of $\Pi$ or,
equivalently, by the $(n+1)$-simplices of the triangulation
of $\Delta$.
\end{defn}

\begin{prop}
\label{prim-uj}
For each $U_j$ there exists $M_j\in ASL_{n+1}(\Z)$ such that
$M_j(U_j)\subset\Sigma_n$ is an open set in the primitive
complex $\Sigma_n$ from Example \ref{si}.
\end{prop}
\begin{proof}
This proposition also follows from the duality with
a unimodular triangulation $\D$ of $\Delta$.
Let $U_j$ be a primitive piece. It corresponds to
a simplex of volume $\frac{1}{(n+1)!}$ in $\D$.
There is an element of $SL_{n+1}(\Z)$ which takes
this simplex to the standard simplex $\Delta_1^{n+1}$
(see \eqref{sta}). Then the image of $U_j$ by the adjoint
to the inverse of this element
is contained in a dual $\Delta_1$-complex.
Such a complex is the result of a translation of $\Sigma_n$ by
Proposition \ref{prim}.
\end{proof}

Recall that a polyhedral complex $\Pi$ is called {\em generic}
at a point $x\in\Pi$ of an open $k$-cell if
$x$ has a neighborhood
homeomorphic to $\R^k\times\Sigma_{n-k}$.

Thus, Proposition \ref{genus} implies that a maximal dual $\Delta$-complex
is a generic polyhedron. In topology such polyhedra often appear
as the so-called {\em special spines} of smooth manifolds.
In the next section we see that $\Pi$
can be compactified so as to become a spine of the polyhedron $\Delta$
after puncturing it in the interior lattice points.

\subsection{Toric varieties and compactification of balanced polyhedra}
\label{compactif}
Consider the complex algebraic torus $\tor$, where $\C^*=\C\setminus 0$.
It is a commutative Lie group under multiplication.
The 2-form
\begin{equation}
\label{symtor}
\frac{1}{2i}
\sum\limits_{j=1}^{n+1} \frac{dz}{z}\wedge\frac{d\bar{z}}{\bar{z}}
\end{equation}
is an invariant symplectic form on $\tor$. There is an action of
the real torus $T^{n+1}=S^1\times\dots\times S^1$
on $\tor$ by coordinatewise multiplication (we treat $S^1\subset\C^*$
as the unit circle). The action of $T^{n+1}$ is Hamiltonian and
thus we have a well-defined {\em moment map} (we refer to \cite{At}
for the general definition or to a textbook, e.g. \cite{CS})
$\Log:\tor\to\R^{n+1}$
\begin{equation}
\label{log}
\Log(z_1,\dots,z_{n+1})=(\log|z_1|,\dots,\log|z_{n+1}|).
\end{equation}

Let $\Delta\subset\R^{n+1}$ be a convex polyhedron
with integer (from $\Z^{n+1}$) vertices. Recall (see
e.g. \cite{GKZ}) that there is a complex toric variety
$\C T_\Delta\supset\tor$.
One way to construct it is to consider {\em the Veronese
embedding} $\tor\to\cp^{\#(\Delta\cap\Z^{n+1})-1}$ defined
by the linear system of monomials associated to $\Delta\cap\Z^{n+1}$.
Here we associate to a point $(p_1,\dots,p_{n+1})$ a monomial
$z^{p_1}\dots z_{n+1}^{p_{n+1}}$. We define $\C T_{\Delta}$
as the closure of the image of the Veronese embedding.
Note that the standard, {\em Fubini-Study}, symplectic form
on the ambient space
$\cp^{\#(\Delta\cap\Z^{n+1})-1}$ defines a symplectic form
on $\C T_\Delta$ (as long as the variety $\C T_\Delta$ is
non-singular). In particular, it gives a symplectic
form $\omega_{\Delta}$ on $\tor$ that is invariant with respect to
the action of $T_\Delta$. This gives us a moment map
with respect to $\omega_\Delta$
$$\mu_\Delta:\tor\to\Delta,
\mu_\Delta(z)=
\frac{1}{\sum\limits_{j\in\Delta\cap\Z^{n+1}} |z^{2j}|}
\sum\limits_{j\in\Delta\cap\Z^{n+1}} j |z^{2j}|.$$
The image of this embedding is the interior $\Int\Delta$.
The map $\mu\Delta$ can be compactified to the moment map
$\bar\mu_\Delta:\C T_\Delta\to\Delta$.

The maps $\Log:\tor\to\R^{n+1}$ and $\mu_\Delta:\tor\to\Int\Delta$
both have the orbits of $T^{n+1}$ as their fibers.
Thus, they define a natural reparametrization
$$\Phi_\Delta:\R^{n+1}\to\Int\Delta.$$

\begin{defn}
Let $\Pi\subset\R^{n+1}$ be an $n$-dimensional
balanced polyhedral complex. By Proposition \ref{obratno}
there is a convex lattice polyhedron $\Delta$ dual to $\Pi$.
We define $\bar{\Pi}\subset\Delta$,
{\em the compactification of $\Pi$},
by taking the closure of $\Phi_{\Delta}(\Pi)$ in $\Delta$.
We call $\bar\Pi\setminus\Phi_\Delta(\Pi)$ {\em the boundary}
of $\bar{\Pi}$. For convenience from now on we identify
$\Pi$ and $\Phi_{\Delta}(\Pi)$.
\end{defn}

\begin{prop}
\label{grani}
Let $\Pi$ be a dual $\Delta$-complex and let
$\Delta'\subset\Delta$
be a $(k+1)$-dimensional face.
Then the intersection $\bar{\Pi}\cap\Delta'$
is a compactification of a dual $\Delta'$-complex $\Pi'$.
If $\Pi$ is maximal then $\Pi'$ is also maximal.
\end{prop}
We prove this proposition simultaneously with the following
proposition describing the behavior of $\Pi$ near infinity.
Recall that a {\em supporting vector} $\stackrel{\to}{v}$
at a face $\Delta'\subset\Delta$ is a vector such that
$p_{\stackrel{\to}{v}}|_\Delta$ reaches its maximum
precisely over $\Delta'$, where $p_{\stackrel{\to}{v}}$
is the orthogonal projection
in the direction of $\stackrel{\to}{v}$.
\begin{prop}
The complex $\Pi'$ from Proposition \ref{grani}
can be obtained in the following way.
Let $L\subset\R^{n+1}$ be the linear $(k+1)$-subspace
parallel to the face $\Delta'$.
Let $\stackrel{\to}{v}$ be a supporting vector at $\Delta'$.
For a sufficiently large $R>0$ we have
$\Pi'=(\Pi-R\hspace{-5pt}\stackrel{\to}{v})\cap L$.
\end{prop}
\begin{proof}
From the finiteness condition in Definition \ref{prpc}
we have that the complex
$\Pi'=(\Pi-R\hspace{-5pt}\stackrel{\to}{v})\cap L\subset L$
does not depend on the choice of $R>0$ and $\stackrel{\to}{v}$
as long as $\stackrel{\to}{v}$ is supporting and $R$ is sufficiently
large. The proof of Proposition \ref{obratno} ensures that
$\Pi'$ is a dual $\Delta'$-complex.
If $\Pi$ is maximal then it is dual to a triangulation of $\Delta$
into simplices of minimal volume. Such a triangulation induces
a triangulation into simplices of minimal volume on the faces
$\Delta'$ and thus $\Pi'$ is also maximal.
\end{proof}

If $\Pi$ is a maximal dual $\Delta$-complex
then it is generic everywhere
except at the points of its boundary $\dd\Pi$.
The following proposition describes the local topology
of $\bar{\Pi}$ near the boundary. It is a corollary
of Proposition \ref{grani}.

\begin{prop}
\label{strata}
Suppose that $\Pi$ is a maximal dual $\Delta$-complex.
A point $x$ in $\bar\Pi$ has a neighborhood
of one of the following $\frac{(n+1)(n+2)}{2}$ types:
$\R^k\times\Sigma_{l-k}\times [0,+\infty)^{n-l}$, where $k\le l\le n$.
Here $k$ is the dimension of the open cell of $\bar\Pi$
which contains $x$ while $l+1$ is the dimension
of the open face of $\Delta$ which contains $x$.
\end{prop}
We call a point with such a neighborhood {\em a $(k,l)$-point}
of $\bar\Pi$.

\begin{rem}
The concept of generic polyhedron is closely related to that
of {special spine} in Topology.
We remind its definition.
Let $M$ be a compact $(n+1)$-manifold with
boundary and $\bar\Pi\subset M$
be an $n$-dimensional CW-complex such that every open
cell is smoothly embedded to $M$.
The complex $\bar\Pi$ is called {\em a spine} of $M$ if $\bar\Pi$ is
a deformational retract of $M$.
The spine $\bar\Pi$ is called {\em special} if for any point
$x\in\bar\Pi\setminus\dd M$
from an open $k$-cell there exists a neighborhood
isomorphic to $\R^k\times\Sigma^{n-k}$.

Note that if $\Int\Delta\cap\Z^{n+1}=\emptyset$
then all the triangulation vertices of a dual $\Delta$-polyhedron
$\Pi$ are from $\dd\Delta$ then $\bar\Pi$ is a spine of $\Delta$.
In general, $\bar\Pi$ is a spine of
the polyhedron $\Delta$ minus a small neighborhood of the
interior lattice points.
Note that $\bar\Pi$ can be treated as a special
spine of $\Delta$ if we treat $\Delta$
as a manifold with corners.
\end{rem}

\subsection{Stratified fibrations}
Let $V$ and $F$ be smooth manifolds, $\Delta\subset\R^{n+1}$
be a lattice polyhedron of full dimension
and $\Pi$ be a maximal dual $\Delta$-complex.

\begin{defn}
\label{strfib}
A smooth map $\lambda:V\to\bar\Pi$ is called
{\em a stratified $F$-fibration} if
\begin{itemize}
\item
The restriction of $\lambda$ to any open $n$-cell $e\subset\bar\Pi$
is a trivial fibration with the fiber $F$;
\item
for each integer pair $(l,k)$,
$0\le k\le l\le n$ there exists a smooth ``model'' map
$\lambda_{l,k}:V_{l,k}\to\Pi_{l,k}$,
where $\Pi_{l,k}\approx
\R^k\times\Sigma_{l-k}\times [0,+\infty)^{n-l}$, such that
any $(l,k)$-point of $\bar\Pi$
has a neighborhood $U\supset x$ such that
$$\lambda|_U:\lambda^{-1}(U)\to U$$
is diffeomorphic to the model map. The model map depends only
on $l$ and $k$.
\end{itemize}
The map $\lambda_{l,k}$ is called {\em the $(l,k)$-fiber
degeneration}; the fiber $F_{l,k}=\lambda^{-1}_{l,k}(x)$ is called
{\em the $(l,k)$-fiber} of $\lambda$.
\end{defn}


The following proposition is a direct corollary of Definition
\ref{strfib}.
\begin{prop}
Let $\lambda:V\to\bar\Pi$ be a stratified fibration
over the compactification of a maximal dual $\Delta$-complex $\Pi$.
For any open $(l,k)$-cell $e$ of $\bar\Pi$
the restriction of $\lambda$ to $e$ is a trivial
fibration over $e$ with the fiber $F_{l,k}$.
\end{prop}

\begin{rem}
Definition \ref{strfib} can be generalized in a straightforward
way to the case when the base is any space with a prescribed
stratification. Here we used the stratification given by
Proposition \ref{strata}.
\end{rem}

\subsection{Hypersurfaces in toric varieties}
Let $f:\tor\to\C$ be a Laurent polynomial
$$f(z)=\sum\limits_j a_jz^j,$$ where $z\in\tor$
and $j\in\Z^{n+1}$ is a multi-index.

We recall that {\em the Newton polyhedron $\Delta$ of $f$}
is the convex hull in $\R^{n+1}$ of the set of all
indices $j\in\Z^{n+1}$ such that $a_j\neq 0$.
Since by assumption $f$ is a polynomial this set is
finite and $\Delta$ is a bounded convex lattice polyhedron.
We also call $\Delta$ the Newton polyhedron of
the hypersurface $\V=\{z\in\tor\ |\ f(z)=0\}$.
According to \cite{GKZ} we call the image $\Log(\V)\subset\R^{n+1}$
{\em the amoeba} of $\V$.

For the rest of the paper we assume that $\Delta$ has a non-empty
interior in $\R^{n+1}$. Otherwise after a suitable
(multiplicative) change of coordinates the polynomial
$f$ can be transformed to a polynomial in a smaller number
of variables.


Let $\C T_{\Delta}$ be the complex toric variety (see e.g. \cite{GKZ})
associated to $\Delta$.
We define $V$ as the closure
of the hypersurface $\V=\{z\in\tor\ |\ f(z)=0\}$ in $\C T_\Delta$.
Taking the Newton polyhedron for $\Delta$ is a canonical choice.
Of course, we can take such compactification
for any convex lattice $(n+1)$-polyhedron $\Delta$,
even if it was not the Newton polyhedron of $\V$.
However the choice of the Newton polyhedron
of $\V$ as $\Delta$ produces the best results as the next
proposition shows.
Recall that in the toric construction there is a $k$-dimensional
complex toric subvariety $\C T_{\Delta'}$ associated to any
$k$-dimensional face $\Delta'\subset\Delta$.

\begin{prop}
\label{NPbest}
The hypersurface $V$ is disjoint from the points
(i.e. the 0-dimensional toric varieties) corresponding to the
vertices of $\Delta$, but intersects all the tori corresponding
to any positive-dimensional face of $\Delta$.

Furthermore, this property characterizes $\C T_{\Delta}$
in the following sense.
Let $\bar\Delta$ be a convex lattice polyhedron in $\R^{n+1}$
with a non-empty interior and $\bar{V}$ be the closure of $\V$ in
$\C T_{\bar{\Delta}}\supset\tor$.
If a hypersurface $\bar{V}$ is disjoint
from the points corresponding to the
vertices of $\bar{\Delta}$ but intersects all the tori corresponding
to positive-dimensional faces of $\bar{\Delta}$
then $\C T_{\bar{\Delta}}=\C T_{\Delta}$.
\end{prop}
\begin{rmk}
Note that even though $\C T_{\Delta}$ is unique
by this proposition, the polyhedron $\bar{\Delta}$ itself
is not unique even up to a translation. The image of $\Delta$
by a homothety with an integer coefficient for $\bar{\Delta}$
corresponds to the same toric variety.
\end{rmk}
\begin{proof}
Proposition \ref{NPbest} follows from the following Lemma.
Note that a vertex is a 0-face of $\Delta$.
\end{proof}
\begin{lem}
Let $\Delta'\subset\Delta$ be a face.
The intersection $V\cap\C T_{\Delta'}$ coincides
with the hypersurface cut on $\C T_{\Delta'}$ by
the closure of the zero set of the following $\Delta'$-truncation
of the polynomial $f$
$$f_{\Delta'}(z)=\sum\limits_{j\in\Delta'}a_j z^j.$$
\end{lem}
\begin{proof}
To prove the lemma it suffices to note that
the monomials from $\Z^{n+1}\cap\Delta'$ have
higher order of vanishing when $z\to\C\Delta'$.
\end{proof}
\begin{rem}
The property of $V$ from Proposition \ref{NPbest}
can be alternatively reformulated in terms of the
moment map $\bar\mu_\Delta:\C T_{\Delta}\to\Delta$,
see subsection \ref{compactif}.
{\em The image $\mu (V)$ is disjoint from the vertices
of $\Delta$ but intersects every positive-dimensional
face of $\Delta$.} According to \cite{GKZ} the image
$\mu(V)$ is called {\em the compactified amoeba} of $\V$.
This restatement is equivalent to the property from
Proposition \ref{NPbest}, since for any face $\Delta'\subset\Delta$
we have $\mu(\C T_{\Delta'})=\Delta'$.
\end{rem}

\begin{figure}[h]
\centerline{\psfig{figure=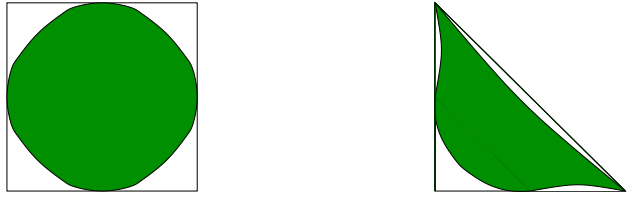,height=1in,width=3.25in}}
\caption{\label{NP} Images of the hyperbola $zw+z+w-1=0$
under the moment maps corresponding to its Newton polygon
and another polygon.}
\end{figure}

\begin{exa}
Let $f(z,w)=zw+z+w-1$. Then $\V\subset(\C^*)^2$ is
a hyperbola. The Newton polygon
$\Delta$ is a square $\{(x,y)\in\R^2\ |\ 0\le x\le 1,0\le y \le 1\}$
and the corresponding toric
surface $\C T_{\Delta}$ is the hyperboloid
$\cp^1\times\cp^1$.


Take now $\bar{\Delta}=\{(x,y)\in\R^2\ |\ 0\le x,0\le y,x+y\le 1\}$.
The corresponding toric surface is $\cp^2\supset(\C^*)^2$.
The images of $\V$ under the associated moment maps
are sketched on Figure \ref{NP}.
\end{exa}

The following example treats projective hypersurfaces.
\begin{exa}
\label{cp}
Let $V\subset\cp^{n+1}\supset\tor$ be a projective hypersurface
of degree $d$ not passing through the points $[1:0:\dots:0],\dots,
[0:\dots:0:1]$.
Then $\V=V\cap\tor$ is given by a polynomial $f$ whose
Newton polyhedron is
$$\Delta_d=\{(x_1,\dots,x_{n+1})\in\R^{n+1}\ |\
0\le x_j, \sum\limits_j x_j\le d\}.$$
Vice versa, $\C T_{\Delta}=\cp^{n+1}$
and the closure of $\V$ in $\cp^{n+1}$ is $V$.
\end{exa}

\subsection{Pairs-of-pants in higher dimensions}

\begin{defn}
\label{pp}
Let $\HH\subset\cp^n$ be the union of $n+2$ generic hyperplanes
in $\cp^{n}$.
Let $\U\subset\cp^n$ be the union of
their $\epsilon$-neighborhoods
for a very small $\epsilon>0$.

The complement $\bar{\PP}_n=\cp^n\setminus\U$ is a manifold with corners.
We call $\bar\PP_n$ {\em the $n$-dimensional pair-of-pants}.
We call $\PP_n=\cp^n\setminus\HH$ {\em the $n$-dimensional
open pair-of-pants}
\end{defn}
Immediately we have the following proposition.
\begin{prop}
A pair-of-pants is a compact manifold with boundary.
An open pair-of-pants is diffeomorphic to the pair-of-pants
minus its boundary.
\end{prop}

\begin{rmk}
Note that the choice of $n+2$ generic hyperplane in
$\cp^n$ is unique up to the action of $PSL_{n+1}(\C)$.
Thus $\PP_n$ can be given a canonical complex structure.
\end{rmk}

Note that $\PP_1$ is diffeomorphic to the Riemann sphere
punctured 3 times, while $\bar\PP_1$ is diffeomorphic
to a closed disk with 2 holes. Thus Definition \ref{pp}
agrees with the classical, one-dimensional, pair-of-pants
definition.

The following proposition describes a natural stratification
of the boundary $\dd\bar\PP$.
\begin{prop}
\label{ddpp}
We have the following canonical decomposition of the boundary
$\dd\bar{\PP_n}=\bigcup\limits_{j=0}^{n-1}\dd_j\bar{\PP_n}$,
where $\dd_j\bar{\PP_n}$ is a $(2n-j)$-dimensional smooth
manifold such that each one of its connected components
is a trivial $T^j$-fibration over $\PP_{n-j}$ (recall that $T^j$
is a $j$-dimensional torus $S^1\times\dots\times S^1$).
Different parts do not intersect:
$\dd_j\bar{\PP_n}\cap\dd_k\bar{\PP_n}=\emptyset$, if $j\neq k$,
but the closure of $\dd_j\bar{\PP_n}$ contains $\dd_k\bar{\PP_n}$
for all $k\le j$.
The number of connected components of $\dd\bar{\PP_n}$
is $\begin{pmatrix} n+2 \\ j+2   \end{pmatrix}$.
\end{prop}
\begin{proof}
Connected components of the manifold $\dd_j\bar{\PP_n}$
can be obtained as the intersections of the boundaries
of the $\epsilon$-neighborhoods of $j$ different
hyperplanes from $\HH$.
\end{proof}


\section{Statement of the results}
\label{mainres}
Let $V\subset\cp^{n+1}$ be a smooth hypersurface of degree $d$.
We choose homogeneous coordinates $[Z_0:\dots:Z_{n+1}]$
so that $V$ is transverse to coordinate hyperplanes
$Z_j=0$ and all their intersections.
The complement of the coordinate hyperplanes in $\cp^{n+1}$
is $\tor$.
Denote $\V=V\cap\tor$.
Then the hypersurface $\V\in\tor$
is given by the equation $f(z)=0$,
where $$z=(z_1,\dots,z_{n+1})=(Z_1/Z_0,\dots,Z_{n+1}/Z_0)$$
stands for affine coordinates in $\tor$ and $f$ is a polynomial
with the Newton polyhedron $\Delta_d$ from Example \ref{cp}.
Recall that we denote the real $n$-dimensional torus
with $T^n=S^1\times\dots\times S^1$.

\begin{thm}
\label{diff}
For every maximal dual $\Delta_d$-complex $\Pi$
there exists a stratified $T^n$-fibration
$\lambda:V\to\bar\Pi$ .
This fibration satisfies to the following properties
\begin{itemize}

\item the induced map $\lambda^*: H^n(\bar\Pi;\Z)
\to H^n(V;\Z)$ is injective,
where $H^n(\bar\Pi;\Z)\approx\Z^{p_g}$, $p_g=h^{n,0}$ is the
geometric genus of $V$;
\item for each primitive piece $U_j$ of $\Pi$ (see Definition \ref{uj})
the inverse image $\lambda^{-1}(U_j)$ is an open pair-of-pants $\PP_n$.
\item for each $n$-cell $e$ of $\bar\Pi$ there exists
a point $x\in e$ such that the fiber $\lambda^{-1}(x)$ is
a Lagrangian $n$-torus $T^n\subset V$;
\item there exist Lagrangian embedding
$\phi_k:S^n\to V$, $k=1,\dots,p_g$
such that the cycles $\lambda(\phi_k(S^n))$ form a basis of $H_n(\bar\Pi)$.
\end{itemize}
Maximal dual $\Delta_d$-complexes exist for every degree $d$ and
every dimension $n$.
\end{thm}

\begin{cor}
\label{lagr}
A $2h^{n,0}$-dimensional subspace of $H_n(V)$ has a basis
represented by embedded Lagrangian tori and spheres.
\end{cor}

Theorem \ref{diff} admits a straightforward generalization to toric
varieties other than $\cp^{n+1}$.
Let $\Delta$ be a bounded convex lattice polyhedron
such that all singularities of the toric variety $\C T_{\Delta}$
are isolated.
Note that the isolated singular points of $\C T_{\Delta}$
necessarily correspond to some vertices of $\Delta$.
Consider the space $(\C^*)^{\#(\Delta\cap\Z^{n+1})}$
of all polynomials of the type
$f(z)=\sum\limits_{j\in\Delta}a_jz^j$ such that $a_j\neq 0$.
Then for a generic choice of a polynomial $f$ from this
space the closure $V$ in $\C T_\Delta$ of the zero set of $f$
is a smooth hypersurface transverse to all toric subvarieties
$\C T_{\Delta'}$ corresponding to the faces $\Delta'\subset\Delta$.
All such $V$ are diffeomorphic and, if we equip them with
the symplectic form from $\C T_{\Delta}$, are symplectomorphic
varieties.

\renewcommand{\thethm}{{1'}}
\begin{thm}
\label{diffg}
For every maximal dual $\Delta$-complex $\Pi$
there exists a stratified $T^n$-fibration
$\lambda:V\to\bar\Pi$.
This fibration satisfies to the following properties
\begin{itemize}
\item the induced map $\lambda^*: H^n(\bar\Pi;\Z)
\to H^n(V;\Z)$ is injective,
where $H^n(\bar\Pi;\Z)\approx\Z^{p_g}$, $p_g=h^{n,0}$ is the
geometric genus of $V$;
\item for each primitive piece $U_j$ of $\Pi$ (see Definition \ref{uj})
the inverse image $\lambda^{-1}(U_j)$ is an open pair-of-pants $\PP_n$.
\item for each $n$-cell $e$ of $\bar\Pi$ there exists
a point $x\in e$ such that the fiber $\lambda^{-1}(x)$ is
a Lagrangian $n$-torus $T^n\subset V$;
\item there exist Lagrangian embeddings
$\phi_k:S^n\to V$, $k=1,\dots,p_g$
such that the cycles $\lambda(\phi_k(S^n))$ form a basis of $H_n(\bar\Pi)$.
\end{itemize}
\end{thm}
\renewcommand{\thethm}{2}

By Remark \ref{nemax} not all convex lattice polyhedra have
maximal dual complexes.
However, in the case of $\Delta_d$ (the polyhedra corresponding to
the projective space), such subdivisions exists for any $d$.
Maximal subdivisions also exist for products of different $\Delta_d$
(this corresponds to hypersurfaces in the product of projective
spaces). It is conjectured that for any lattice polyhedron $\Delta$
there exists a sufficiently large integer $N$ that $N\Delta$
(the result of scaling of $\Delta$ by $N$)
has a maximal subdivision.

The next theorem describes the behavior of the fibration
$\lambda$ with respect to a complex structure
on $V$. Recall that, unlike the smooth and symplectic structures,
the complex structure on $V$ depends on the polynomial
$f$ and not just on $\Delta$.

Recall that a map $\lambda:V\to\bar\Pi$
is called {\em a totally real fibration} if for any
$z\in V$ the tangent space to the fiber through $z$
is totally real i.e. contains no positive-dimensional complex
subspaces (as long as the fiber is smooth near $z$).
We say that a hypersurface $V\subset\C T_{\Delta}$
is {\em defined over $\R$} if it can be obtained as the
closure of the zero set of a polynomial $f:\tor\to\C$
whose coefficients are real.

\begin{thm}
\label{comp}
For every maximal dual $\Delta$-complex $\Pi$
there exists a smooth hypersurface $V\subset\C T_{\Delta}$
defined over $\R$ such that the map $\lambda$ from
Theorem 1' preserves the real structure  of $V$,
i.e. $\lambda\circ\conj=\lambda$, where $\conj:V\to V$ is the
involution of complex conjugation.
Furthermore, $\lambda$ is a totally real fibration.
\end{thm}
\renewcommand{\thethm}{3}

Theorems 1' and \ref{comp} can be extended further
to polyhedra $\Delta$ corresponding to toric varieties with non-isolated singularities.
However, in order to do that,
one has to modify the definition of stratified fibrations
to include singular total spaces $V$.
We do not do that. In the next theorem
we no longer have any restrictions on the
convex lattice polyhedron $\Delta$, but
its statement concerns only the toric, non-singular,
part $\V\subset\tor$ of the hypersurface $V$.

\begin{thm}
\label{tor}
For every maximal dual $\Delta$-complex $\Pi$
there exists a stratified $T^n$-fibration
$\lambda^\circ:\V\to \Pi$.
This fibration satisfies to the following properties
\begin{itemize}
\item the induced map $(\lambda^\circ)^*: H^n(\Pi;\Z)
\to H^n(\V;\Z)$ is injective,
where $H^n(\Pi;\Z)\approx\Z^{p_g}$, $p_g=h^{n,0}$ is the
geometric genus of $V$;
\item for each primitive piece $U_j$ of $\Pi$ (see Definition \ref{uj})
the inverse image $(\lambda^\circ)^{-1}(U_j)$ is an open pair-of-pants
$\PP_n$.
\item for each $n$-cell $e$ of $\Pi$ there exists
a point $x\in e$ such that the fiber $(\lambda^\circ)^{-1}(x)$ is
a Lagrangian $n$-torus $T^n\subset V$;
\item there exist Lagrangian embeddings
$\phi_k:S^n\to\V$, $k=1,\dots,p_g$
such that the cycles $\lambda^\circ(\phi_k(S^n))$ form a basis of $H_n(\Pi)$.
\end{itemize}
\end{thm}

          \
\begin{rmk}
These theorems generalize to complete intersections.
The base of the fibration in this case is the intersection
of the maximal dual balanced polyhedra for the corresponding
hypersurfaces (we have to choose them in a mutually
general position).

From a different point of view the base is dual to
a maximal mixed lattice subdivision of
the Newton polyhedra of the participating equations.
The primitive pieces for complete intersections are
products of the primitive pieces for hypersurfaces.
Sturmfels' generalization \cite{St} of the
patchworking technique allows to produce in this case
the Lagrangian lifts
of the base cycles.

This generalization will be the subject of a future paper.
\end{rmk}

\section{Some examples}

\subsection{Riemann surfaces}
\label{RS}
Let $S$ be a closed Riemann surface of genus $g>1$.
It is well-known that $S$ admits a decomposition into
pairs-of-pants.
Namely, there exist $3g-3$ disjoint embedded circles
$C_j\subset S$ such that
$S\setminus \bigcup\limits_{j=1}^{3g-3} C_j$ is
a disjoint union of $2g-2$ copies of the pair-of-pants $P$.
The pair-of-pants surface $P$ is homeomorphic to
the Riemann sphere $\cp^1$ punctured in three points.

To such a decomposition we associate a graph $\Gamma$.
The vertices of $\Gamma$ correspond to the pairs-of-pants
while the edges correspond to the circles $C_j$.
Each edge joins the vertices corresponding to the
adjacent pairs-of-pants.

\begin{figure}[h]
\centerline{\psfig{figure=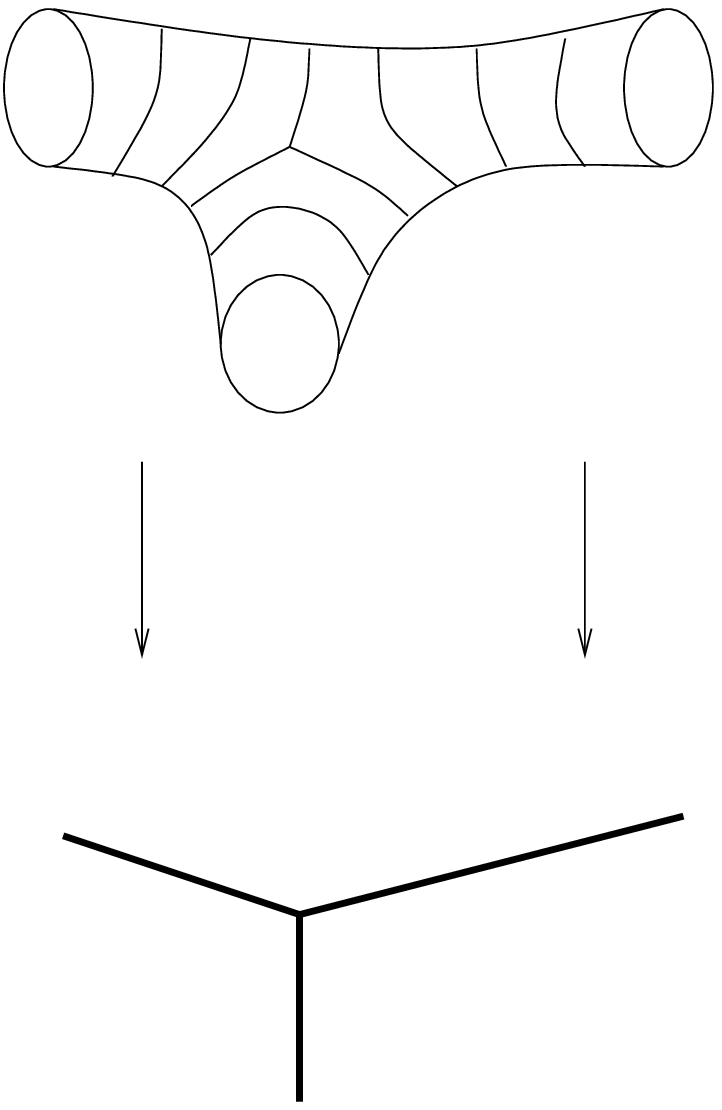,height=2.2in,width=1.5in}
\hspace{0.5in}
\psfig{figure=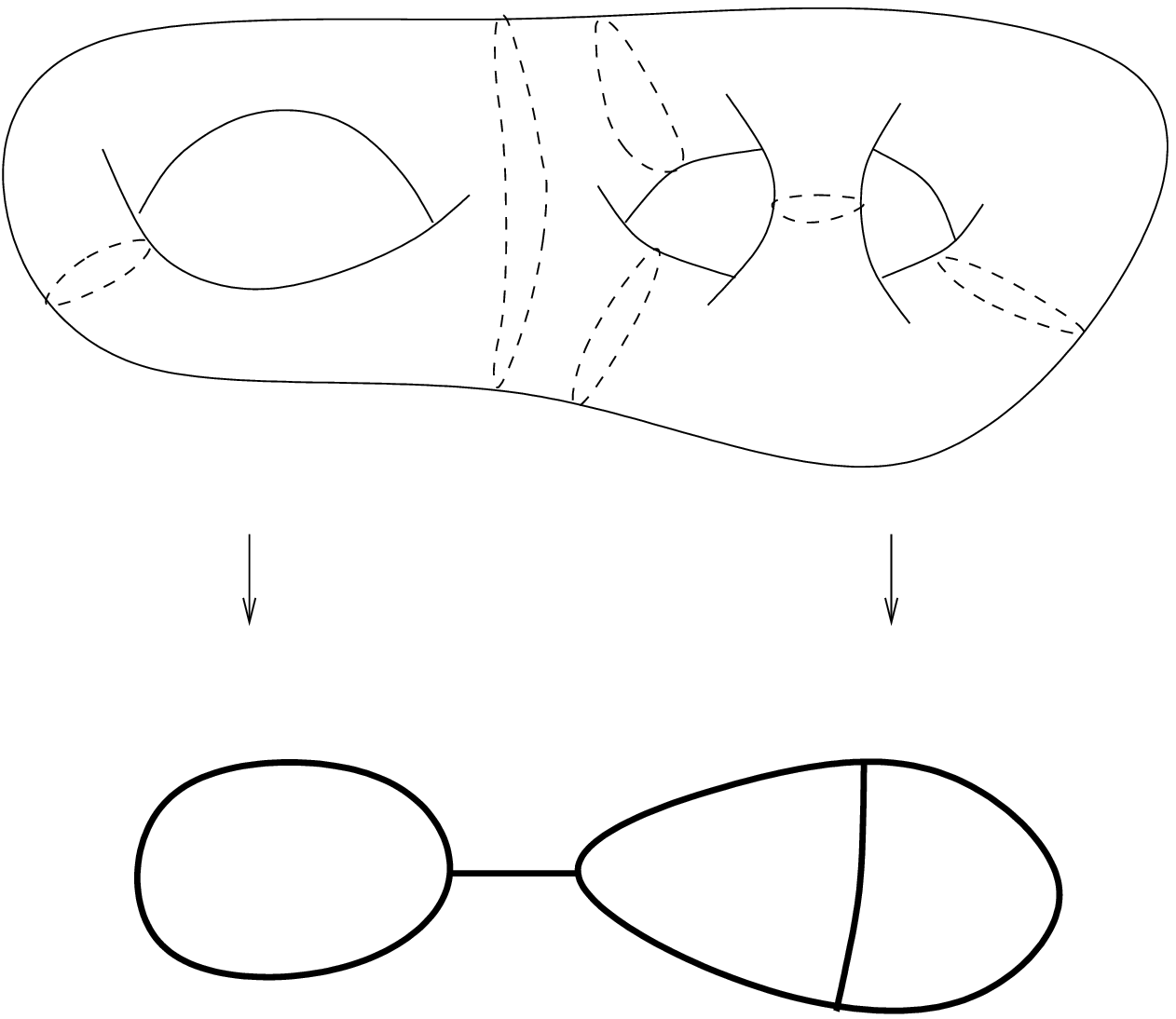,height=2.2in,width=2.5in}}
\caption{\label{tripod} Circle fibrations on a pair-of-pants
and on a surface with a pair-of-pants decomposition.}
\end{figure}

There exists a fibration $\pi:S\to\Gamma$ such that
the circles $C_j$ are inverse images of the midpoints
of the edges of $\Gamma$. Such fibration is canonically
associated to our decomposition into pairs-of-pants.
To construct it we fiber each individual pair-of-pants
over a tripod graph as pictured on the left-hand-side
of Figure \ref{tripod}. Corresponding diagrams in the Newton
polygon of a polynomial were explored in \cite{RuN}.

\subsection{The elliptic curve and the K3-surface}
Here we consider the well-known fibrations of the elliptic
curve and the K3-surface.

Let $\C E$ be an elliptic curve, i.e. a Riemann surface
of genus 1. Since $\C E$ is topologically a torus,
there is a trivial $S^1$-fibration $\lambda_E:\C E\to S^1$.

Suppose that the elliptic curve $\C E\subset\C T_\Delta$ is
presented as a curve in a toric surface $\C T_\Delta$,
where $\Delta$ is the Newton polygon of a polynomial defining $\C E$.
By the genus formula (see \cite{Kh}), $\Int\Delta$ contains a
unique lattice point. By Proposition \ref{genus}
a dual $\Delta$-complex is homotopy equivalent to a circle.
It is easy to see that the fibration from Theorem \ref{diffg}
coincides up to homotopy with the trivial $S^1$-fibration
$\C E\approx S^1\times S^1\to S^1$.

Another famous fibration $\lambda_K:\C K\to S^2$
has the K3-surface $\C K$ as its total space.
All its fibers, except for 24 of them are Lagrangian tori.

Suppose that the polyhedron $\Delta$ has exactly one interior
lattice point. Then, by Khovanskii's formula \cite{Kh},
the zero locus $\C K$ of a generic polynomial with
the Newton polyhedron $\Delta$ is a K3-surface.
A dual $\Delta$-complex is homotopy equivalent to a sphere $S^2$
by Proposition \ref{genus}.

Again, the fibration $\lambda$ can be deformed to a fibration
like $\lambda_K$ by so-called shelling of $\bar\Pi$
\footnote{A higher-dimensional version of such deformation
will be the subject for a future paper.}.

In higher dimensions, if $\Delta$ is a non-singular
polyhedron with a unique interior
lattice point, then the corresponding hypersurface
$V\subset\C T_\Delta$ is a smooth Calabi-Yau manifold.
Singular torus fibrations $V\to S^n$ were constructed
by Zharkov \cite{Zh}.
Ruan \cite{Ru} noted that such fibrations can be made
Lagrangian.

Theorem \ref{diffg} constructs in this case
a stratified torus fibration over a polyhedral complex
homotopy equivalent to $S^n$.

\subsection{Hyperplanes in the projective space}
\label{hyperplanes}
This is a fundamental example for the main theorems.
Let $H=\{z_1+\dots+z_{n+1}+1=0\}\subset\cp^{n+1}$
be a hyperplane. Its toric part
$H^{\circ}=H\cap\tor$ is an open pair-of-pants.

Let $\Log$ be the moment map for $\tor$ (see \eqref{log}).
\begin{lem}
$\Sigma_n\subset\Log(H^{\circ})$.
\end{lem}
\begin{proof}
By \cite{PaRu} $\Sigma_n$ is a spine of the amoeba $\Log(H^{\circ})$
and, therefore, its subset. The lemma can alternatively
be verified by writing explicit inequalities defining $\Log(H^\circ)$.
\end{proof}

The complement $\R^{n+1}\setminus\Sigma_n$ consists of
$n+2$ components. Each component is the region where
one of the functions $0,x_1,\dots,x_{n+1}$ is maximal.
In the component corresponding to $x_j$ we consider the foliation
into straight lines parallel to the gradient of $x_j$
(the $j$th basis vector). In the component corresponding to $0$
we consider the foliation into straight lines parallel to
$(1,\dots,1)$. These foliations glue to a singular foliation
$\F'$ which has singularities at $\Sigma_n$.

\begin{figure}[h]
\centerline{
\psfig{figure=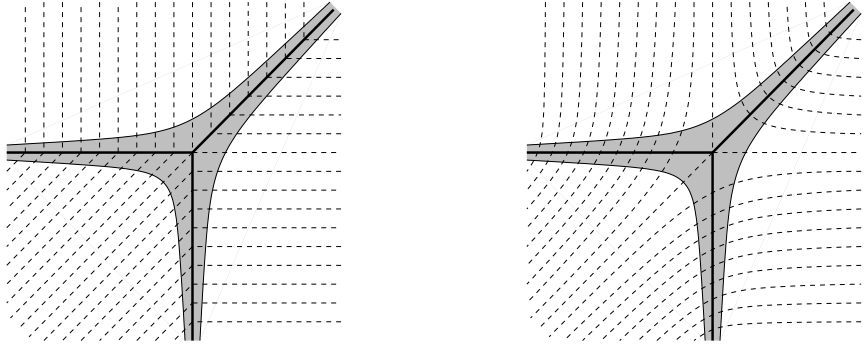,height=1.8in,width=4.5in}}
\caption{\label{fol} The amoeba $\Log(H^\circ)$ together
with the foliation $\F'$ and its deformation $\F$.}
\end{figure}

It is easy to smooth out $\F'$ (in a symmetric way
with respect to the homogeneous coordinates permutations)
at the open $n$-cells of $\Sigma_n$
(see Figure \ref{fol}). However, the singularities at the
smaller-dimensional cells are essential. The leaves passing through
an open $(n-k)$-cell are homeomorphic to the cone over $k+2$ points.


We denote the resulting
foliation with $\F$. The foliation $\F$ is a singular fibration
and defines the projection
$\pi_\F:\R^{n+1}\to\Sigma_n$.

The following statement is a key lemma in the proof of the main
theorems of this paper.
\begin{lem}
\label{thm3hyp}
The composition
$$\lambda_H=\pi_\F\circ\Log:H^\circ\to\Sigma_n$$
is a stratified $T^n$-fibration in the sense of Definition \ref{strfib}.
It satisfies to all conclusions of Theorem \ref{tor} except for
the third one. The fibration $\lambda_H$ can be deformed so that
the third condition will also hold.
\end{lem}
The proof of this lemma occupies the rest of this subsection.

To figure out the fibers of $\lambda_H$ we need to understand
the critical points of $\Log|_{H^\circ}$. Following \cite{Ka-gauss}
and \cite{Mi} for a hypersurface $V^\circ\subset\tor$ we define
{\em the logarithmic Gauss map}
$$\gamma:V^\circ\to\cp^n$$ by
taking the composition of a branch of a holomorphic logarithm
of each coordinate with the conventional Gauss map. This produces
the following formula
$$\gamma(z_1,\dots,z_{n+1})=
[z_1\frac{\dd f}{\dd z_1}:\dots:z_{n+1}\frac{\dd f}{\dd z_{n+1}}],$$
where $f$ is the polynomial defining $V^\circ$.

Note that the Newton polyhedron of $z_j\frac{\dd f}{\dd z_j}$
coincides with the Newton polyhedron $\Delta$ of $f$.
Therefore, by Kouchnirenko's formula \cite{Kou},
$\deg\gamma=(n+1)!\Vol\Delta$. In particular, if $V^\circ=H^\circ$
then $\deg\gamma=1$.

\begin{lem}[cf. Lemma 3 of \cite{Mi}]
The set of critical points of $\Log|_{V^\circ}$ coincides
with $\gamma^{-1}(\R P^n)$.
\end{lem}
\begin{proof}
Let $z\in V^\circ$ and let $\LL$ be a branch of a holomorphic
logarithm $(z_1,\dots,z_{n+1})\mapsto (\log(z_1),\dots,\log(z_{n+1}))$
defined in a neighborhood of $z$. The point $z$ is critical
for $\Log|_{V^\circ}$ iff $V^\circ$ and the orbit of the real
torus $T^n$ are not transversal at $z$. But $\LL$ takes the
tangent space to an orbit of $T^n$ to
a translate of $i\R^{n+1}$ in $\C^{n+1}$.

Therefore, $z$ is critical iff $\LL(T_zV^\circ)$ contains at least
$n$ purely imaginary vectors which is, in turn, equivalent to
$\gamma(z)\in\R P^n$.
\end{proof}

\begin{cor}
The set of critical points of $\Log|_{H^\circ}$ coincides
with the real locus $\R H^\circ$ of $H^\circ$
(i.e. with the set of real solutions of $z_1+\dots+z_{n+1}+1=0$).
\end{cor}
\begin{proof}
Note that, since $H^\circ$ is defined over $\R$, we have
$\gamma(\R H^\circ)\subset\rp^n$. Note that $\gamma$ extends to
a map $H\to\cp^n$ which is an isomorphism, since $\deg\gamma=1$.
\end{proof}


\begin{cor}
The locus $\DD\subset\Log(H^\circ)$ of critical
values of $\Log|_{H^\circ}$
is an immersed manifold transverse to the foliation $\F$.
\end{cor}
\begin{proof}
The map $\Log|_{\R H^\circ}:\R H^\circ\to\DD\subset\R^{n+1}$
is an immersion since the map $\Log|_{\rtor}:\rtor\to\R^{n+1}$
is an immersion (it is a trivial $2^{n+1}$-covering of $\R^{n+1}$).

To see the transversality we recall the definition of
the foliation $\F'$.
For each component of $\R^{n+1}\setminus\Sigma_n$ the foliation
$\F'$ is parallel to a vector $\ve$ normal to a facet $\Delta'$
of the Newton polyhedron of $H^\circ$. Therefore, any hyperplane tangent to
$\gamma(\R H^\circ)$ is transverse to $\ve$
(hyperplanes parallel to $\ve$ corrspond to the intersection
of $\R H$ with the divisor corresponding to $\Delta'$).
Furthermore, hyperplanes close to being parallel to $\ve$ are
close to the hyperplane in $\cp^{n+1}$ corresponding to this facet
and therefore are far from the given component of
$\R^{n+1}\setminus\Sigma_n$. Thus the result $\F$ of smoothing
is also transverse to $\DD$ and the angle between them in $\R^{n+1}$
is separated from 0.
\end{proof}

Note that $\pi_\F$ is a stratified $[-1,1]$-fibration.
Thus, the transversality of $\DD$ and $\F$ implies that $\lambda_H$
is a stratified fibration for $\Sigma_n$.
We need to show that the restriction
of $\pi_F$ to open $n$-cells of $\Sigma$ is a torus fibration.

Consider a point $x=(-t,\dots,-t,0)$ for a large $t>0$.
Note that $\DD$ is almost horizontal near $x$.
Thus the fiber of $\lambda_H$ over $x$ is diffeomorphic
to the fiber $F$ of a composition of $\Log|_{H^\circ}$ and
the linear projection onto the first $n$ coordinates.
Note that the map $F\to T^n$
obtained by taking the arguments of the first $n$ coordinates
is a diffeomorphism.
Recall that $H^\circ$ is given by the equation $z_1+\dots+z_{n+1}+1=0$.
The absolute values of the coordinates $z_1,\dots,z_n$ are fixed.
For any value of their argument we take $z_{n+1}=1-z_1-\dots-z_n$
to get the unique point from $F$ corresponding to this choice
of the arguments. Since $|z_1|,\dots,|z_n|$ are small $z_{n+1}\neq 0$.

We verify the conclusions of Theorem \ref{tor} item-by-item.
The first and the last conclusions are vacuous in this case,
since $\Sigma_n$ (and, therefore,
$\bar\Sigma_n$ as well) is contractible.
The second one holds since $H^\circ$ is itself an open
pair-of-pants.

To make the third conclusion true we have to modify $\lambda_H$
a little. The fiber $F$ is not Lagrangian, but it is close
to a Lagrangian torus $\Lambda=\{|z_j|=\const,\ j=1,\dots,n,\
z_{n+1}=-1\}$. We can deform $H^\circ$ a little in a neighborhood
of $\Log^{-1}(x)$ to make it intersect the fiber of $\pi_{\F}\circ\Log$
along $\Lambda$. Therefore, $F$ is Lagrangian for a nearby
symplectic structure. By Moser's trick (see e.g. \cite{CS})
there exists a self-diffeomorphism $h$ of $H^{\circ}$ constant
outside of a neighborhood of $\Log^{-1}(x)$ and taking
one symplectic structure to another.
We redefine $\lambda$ as $\lambda\circ h$.
This ensures a Lagrangian fiber over one of the
$\begin{pmatrix} n+2 \\ 2 \end{pmatrix}$ open $n$-cells
of $\Sigma_n$. We do the same for all other $n$-cells.

\subsection{A localization $Q^n\subset\tor$ of the standard hyperplane}
\label{lochyp}
The toric part $H^\circ\subset\tor$ of a hyperplane from
\ref{hyperplanes}
is a nice embedding of $\PP_n$ to $\tor$. However
for our purposes it is convenient to modify it in
a neighborhood of infinity to get a different submanifold $Q^n$ which
is better suited for gluing.

Note that the symmetric group $S_{n+2}$
acts on $\cp^{n+1}$ by interchanging the $n+2$
homogeneous coordinates.
This action leaves $\tor$ and $H^\circ\subset\tor$
invariant.
\begin{prop}
\label{locprop}
There exists a proper submanifold $Q^n\subset\tor$
such that
\begin{itemize}
\item $Q^n$ is embedded in $\tor$ symplectically, i.e. so that the
restriction of the form \eqref{symtor} to $Q^n$ is a symplectic form.
\item $Q^n$ is isotopic to $H^\circ$ in $\tor$.
\item The composition $\pi_{\F}\circ\Log_t|_{Q^n}$
is a stratified $T^n$-fibration that satisfies to all
hypotheses of Theorem \ref{tor}.
\item the closure $\bar{Q^n}$ of $Q^n$ in $\cp^{n+1}\supset\tor$
is a smooth manifold isotopic to $H$.
\item $Q^n$ is invariant with respect to the action of the
symmetric group $S_{n+2}$ on $\tor$ (see above).
\item For a sufficiently large $M>0$
$$Q^n\cap\tor_{-M}=
Q^{n-1}\times\\C^*_{-M},$$
where $\tor_{-M}=\{(z_1,\dots,z_{n+1})\in\tor\ |\ \log|z_{n+1}|<-M\}$
and $\C^*_{-M}=\{z\in\C^*\ |\ \log|z|<-M\}$.
In particular, the intersection
$Q^n\cap\tor_{-M}$
is invariant under a translation $z_{n+1}\mapsto cz_{n+1}$,
$0<c<1$.
\end{itemize}
\end{prop}

\begin{figure}[h]
\centerline{\psfig{figure=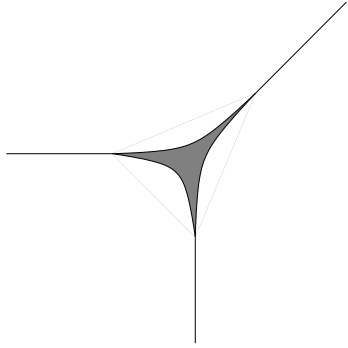,height=1.8in,width=1.8in}}
\caption{\label{loc} The amoeba of the localization $Q^n$
of a hyperplane.}
\end{figure}

\begin{proof}
We construct $Q^n$ inductively by dimension $n$.
If $n=0$ then $H^\circ$ is a point and $Q^0=H^\circ$.
Assume that $Q^k$, $k<n$ is already constructed.
Consider the simplex
$$\Delta_n(R)=\{x\in\R^{n+1}\ |\ -x_j\le R,\sum\limits_j x_j\le R\}.$$
Each its $k$-dimensional face is dual to a $(n+1-k)$-cell of
$\Sigma_n$.
Fix a sufficiently large number $R_n>0$.

First we define $Q^n\cap\Log^{-1}(\dd\Delta(R_n))$.
Each $k$-face of $\Delta(R_n)$ is contained in a unique
affine $k$-space $A$ in $\R^{n+1}$. Furthermore, the adjoint
faces cut the polyhedron $\Delta_{k-1}(R_n)\subset A$.
Thus we may identify $A$ with $\R^k$ and, therefore,
$\Log^{-1}(A)$ with $(\C^*)^k$.
By the induction assumption we already have $Q^{k-1}\subset(\C^*)^k\to\R^k$.
We define $Q^n\cap\Log^{-1}(\dd\Delta(R_n))$ to be equal
to the union of these $Q^k$ for all faces
of $\dd\Delta(R_n)$. By the induction hypothesis (and since
$R_n$ was large enough) the choices
over different faces agree.

Our next step is to extend $Q^n$ to the complement of
$\Log^{-1}(\Delta(R_n))$. For each face $\Delta'$ of
$\dd\Delta(R_n)$ consider its outer normal cone
$C_{\Delta'}\subset\R^{n+1}$ (e.g. if $\Delta'$ is a facet
then $C_{\Delta'}$ is a ray). We define
$$Q^n\cap\Log^{-1}(\Delta'+C_{\Delta'})=
\bigcup\limits_{\ve\in C_{\Delta'}}
e^{\ve} Q^n\cap\Log^{-1}(\Delta').$$
In other words, we span the region above the normal cone
of a $k$-face $\Delta'$ by the translates of the manifold $Q^k$.

We set $Q^n\cap\Log^{-1}(\Delta(R_n-1))=
H^\circ\cap\Log^{-1}(\Delta(R_n-1))$. By now we have defined
$Q^n$ everywhere, but $\Log^{-1}(\Delta(R_n)\setminus
\Delta(R_n-1))$.

Consider a facet $\Delta'$ of $\dd\Delta(R_n-1)$,
e.g. the one sitting in the hyperplane $A=\{x_{n+1}=R_n-1\}$.
Since $R_n$ is large enough, $z_{n+1}^{R_n-1}$ is small
enough and the intersection
$H^\circ\cap\Log^{-1}(A)$ is close enough to
the zero set of $z_1+\dots+z_n+1=0$. By the induction
hypothesis this zero set can be deformed to $Q^{n-1}$.
We define $Q^n\cap\{\Log|z_{n+1}|=t\}$, $-R_n\le t\le -R_n+1$
using this deformation.
We repeat the same procedure for all other facets of
$\Delta(R_n-1)$.
\end{proof}

Denote $\bar{Q}^n=Q^n\cap\Log^{-1}(\Delta(R_n+1))$.
This is the core part of $Q^n$ and is diffeomorphic
to a closed pair-of-pants $\bar{\PP_n}$ (as a manifold
with corners).

\section{Reconstruction of the complex hypersurface from
a balanced polyhedron $\Pi$}
Theorem \ref{diffg} can be treated as a pair-of-pants
decomposition for $V$. We can use this presentation to
reconstruct $V$ from $\Pi$.
This allows to interpret a maximal balanced
polyhedral complex $\Pi$ as the complex
encoding the gluing pattern of pairs-of-pants
in order to get $V$. Here is the way to reconstruct $V$ from $\Pi$.

For each vertex $v_j$ of $\Pi$ take a copy $\bar Q_j$ of $\bar\PP_n$.
This copy can be identified with the localized hyperplane
$\bar Q^n\subset\tor$.
Recall that by Proposition \ref{ddpp} $\dd_1(\bar Q_j)$ consists
of $n+2$ components. Each such component corresponds to a 1-cell
of $\Pi$ adjacent to $v_j$.

Let $e_{jk}$ be a 1-cell of $\Pi$ connecting the vertices
$v_j$ and $v_k$. For each such 1-cell we identify
the closures $F_j$ and $F_k$ of the corresponding components of
$\dd_1(\bar Q_j)$ and $\dd_1(\bar Q_k)$ in the following way.
%

Without the loss of generality we may assume
that in both copies $\bar Q_j, \bar Q_k$ of $\bar Q^n$
the edge $e_{jk}$ corresponds
to the facet $x_{n+1}=-R_n$ of $\Delta(R_n)$ (see \ref{lochyp}).
(Note that such correspondence is given by matrices $M_j,M_k$
from Proposition \ref{prim-uj}.)
We attach $F_j$ to $F_k$ by the map
$$(z_1,\dots,z_n,z_{n+1})\mapsto (z_1,\dots,z_n,\bar{z}_{n+1}),
\log|z_{n+1}|=-R_n,$$
where $\bar{z}_{n+1}$ is the complex conjugate to $z_{n+1}$.

The result $U$ of this gluing is a manifold with boundary.
The boundary comes from the unbounded cells of $\Pi$.
Denote $W^\circ=U\setminus\dd U$.
The boundary is formed by the closures $F$ of the components
of $\dd_1(Q_j)$ that correspond to unbounded 1-cells in $\Pi$.
By Proposition \ref{ddpp} each such $F$ is a circle fibration over
a union of lower-dimensional pairs-of-pants $\PP_{n-1}$.
Let $W$ be the result of collapsing all fibers of these
fibrations on $\dd U$. Note that $W$ is canonically a smooth
manifold since this procedure locally coincides with
collapsing the boundary on $\bar\PP_n$ which results in $\cp^n$.

\renewcommand{\thethm}{4}
\begin{thm}
\label{locthm}
The manifold $W$ is diffeomorphic to $V$.
The manifold $W^\circ$ is diffeomorphic to $\V$.
\end{thm}

\begin{cor}
The manifolds $W$ and $W^\circ$ depend only on the lattice polyhedron
$\Delta$ associated to $\Pi$, not on $\Pi$ itself.
\end{cor}

\begin{rem}
With a little more care this reconstruction process can be
made in the symplectic category, i.e. the result $W$ of gluing
can be given a natural symplectic structure. This is due
to the following two reasons.
The first one is that the pair-of-pants possesses a natural symplectic
structure (the one which gives the standard symplectic $\cp^n$
after the symplectic reduction of the boundary).
The second one is that two pairs-of-pants get identified along
a part $F$ of their boundary which is a symplectically flat
hypersurface, it has a neighborhood $F\times [0,1]$
symplectically isomorphic to $Q_{n-1}\times A$, where $A\subset\C^*$
is an annulus. This product is consistent with the
$S^1$-fibration $F\to Q_{n-1}$ from Proposition \ref{ddpp}.
\end{rem}

\section{Proof of the main theorems}
We are free to choose any smooth hypersurface $V$ with
the Newton polyhedron $\Delta$ to construct the stratified
fibration $\lambda$, since all such hypersurfaces are isotopic.
We use Viro's patchworking construction
\cite{Vi} to choose a convenient $V$. Recall that the Newton
polyhedron $\Delta\subset\R^{n+1}$ of $V$ is a convex
polyhedron whose vertices are lattice points.

\subsection{Viro's patchworking}
Let $v:\Delta\cap\Z^{n+1}\to\R$ be any function and
$a(z)=\sum\limits_{j\in\Delta\cap\Z^{n+1}}a_jz^j$ be any
polynomial.
Following \cite{Vi} we define the
{\em patchworking polynomial} for any $t>0$ by
$$f^v_t(z)=\sum\limits_{j\in\Delta\cap\Z^{n+1}}a_jt^{-v(j)}z^j,$$
where $a_j\neq 0$ for any $j\in\Delta\cap\Z^{n+1}$.
Note that if $v$ is integer-valued then $f^v_t$ makes sense also
for any $t\in\C^*$.

\begin{rem}
\label{patch-orig}
In \cite{Vi} the patchworking polynomial was used for construction
of real algebraic hypersurfaces with controlled topology. The topology
of the zero set of a real patchworking polynomial for $t>>0$ depends only
on the function $v$ and on the signs of the coefficients $a_j$.
\end{rem}

\subsection{Non-Archimedian amoebas}
\label{naa}
If $V\subset\tor$ be an algebraic variety.
The image $\Log(V)\subset\tor$ is called the {\em amoeba} of $V$,
see \cite{GKZ}. Note that amoebas make sense also for varieties
over other fields $K$ as long as we have a norm
$K^*=K\setminus\{0\}\to\R_+$. The map $\Log_K:(K^*)^{n+1}\to\R^{n+1}$
is defined by $\Log_K(z_1,\dots,z_{n+1})=
(\log||z_1||_K,\dots,\log||z_{n+1}||)$ and the amoeba of
$V_K\subset (K^*)^{n+1}$ is defined to be $\Log_K(V_K)$.

A particularly useful case is when $K$ is an algebraically closed
field with a non-Archimedian valuation.
Recall that a non-Archimedian valuation
\footnote{Sometimes a valuation is defined as minus such a function.}
is a function
$\val:K^*\to\R$
such that $\val(a+b)\le\max\{\val(a),\val(b)\}$
and $\val (ab)=\val (a)+\val (b)$. Note that $e^{\val}$
gives a norm on $K$
and $\Log_K$ is nothing but taking the coordinatewise valuation.

Non-Archimedian amoebas of hypersurfaces were completely
described in \cite{Ka}.
An example of such field is the field $K$ of the Puiseux series
with complex coefficients in $t$. Namely an element of $K$
is a formal series $b(t)=\sum\limits_{k\in J}b_kt^k$, $b_k\in\C^*$
where $J\subset\R$ is any bounded from below set contained in
a finite union of arithmetic progressions.
The valuation is defined by $\val ||b(t)||=-\min J$.
Note that we used irrational as well as rational powers
in the Puiseux series to make the valuation surjective.

\begin{thmquote}[Kapranov \cite{Ka}]
If $V_K\subset\ktor$ is a hypersurface given by a
polynomial $f=\sum a_j z^j$, $a_j\in K^*$ then
the (non-Archimedian) amoeba of $V_K$ is the balanced
polyhedral complex corresponding to the function
$v(j)=\val(a_j)$ defined on the lattice points
of the Newton polyhedron $\Delta$ of $V_K$ as in Example \ref{tri}.
\end{thmquote}

\subsection{Lifts of non-Archimedian amoebas to $\tor$}
\label{lifts}
Consider the map $u:K^*\to S^1$ defined by
$u(b)=\arg(b_{-\val(b)})$, $b=\sum\limits_{k\in J}b_kt^k$.
In other words, $u$ takes the argument of the coefficient at the
lowest power of $t$. This is a homomorphism
from the multiplication group $K^*$.
Together with $\val$ it gives a homomorphism
$w=(\val,u):K^*\to\C^*\approx\R\times S^1$ and thus a homomorphism
$W:(K^*)^{n+1}\to\tor$.

\begin{lem}
If $V\subset (K^*)^{n+1}$ is a hypersurface given by a
polynomial $f=\sum a_j z^j$, $a_j\in K^*$ then
$W(V_K)\subset\tor$ depends only on the values $w(a_j)\in\C^*$
of the coefficients.
\end{lem}
\begin{proof}
Kapranov's theorem takes care of $\Log(w(V_K))=\Log_K(V_K)$.
We need to prove that the values $u(a_j)$ determine
the arguments of $W(V_K)$.
Let $x\in\Log_K(V_K)$. By Kapranov's theorem it means
that there is a set of indices $j_1,\dots,j_l$, $l\ge 2$, such that
$\val(a_{j_1})=\dots=\val(a_{j_l})\ge\val(a_j)$ for any
other index $j$. Let $z\in(K^*)^{n+1}$ be a point such that
$\Log_K(z)=x$. The lowest powers of $t$ in the Puiseux
series $f(z)$ are contributed by the monomials $a_{j_1}z^{j_1},
\dots,a_{j_l}z^{j_l}$. If $f(z)=0$ then the coefficients at
these lowest powers are such that their sum is zero. Conversely,
the higher powers of $t$ can be arranged to make $f(z)=0$
without the change of $W(z)$ as in the proof of Kapranov's theorem.
\end{proof}

\subsection{Maslov's dequantization}
Consider the following family of binary operations on $\R\ni x,y$:
$$x\oplus_t y=\log_t(t^x+t^y),$$
for $t>1$ and
$$x\oplus_\infty y=\lim\limits_{t\to 0} x\oplus_t y=\max\{x,y\}.$$
This is a commutative semigroup operation (no inverse
elements and no zero) for each $t$. The set $\R$
equipped with this operation for addition and
with $x\odot y=x+y$ for multiplication is a semiring $\R_t$.
Indeed, for any $x,y,z\in\R$ we have $x\odot(y\oplus_t z)=
(x\odot z)\oplus_t (y\odot_t z)$.

Passing from a finite $t$ to infinity in this family of
semirings is called {\em Maslov's dequantization}, cf. \cite{Ma}.
Note that for all finite values of $t$ the semiring is isomorphic
to the semiring of real positive numbers equipped with the usual
addition and multiplication. But the behavior at $t=\infty$ is
qualitatively different, the addition becomes idempotent,
$x\oplus_\infty x=x$. The prefix ``de" reflects the fact that
in this deformation the classical calculus operations appear
on the quantum side.

There is a universal bound for the convergence of the operations
$\oplus_t$ to $\oplus_\infty=\max$. Namely, we have
\begin{equation}
\label{Md}
\max\{x_1,\dots,x_N\}\le
x_1\oplus_t\dots\oplus_t x_N\le
\max\{x_1,\dots,x_N\}+\log_t N
\end{equation}

The dequantization point of view can be used to reinterpret
Viro's patchworking, see \cite{Vi-deq}. Instead of deforming
the coefficients of the polynomial we may keep them constant,
but deform the addition operation instead. This point of
view yields some useful estimates on the zero set of the
patchworking polynomial as shown below.

One way to think of a polynomial is to think of it
as a collection of coefficients at its monomials.
Fix a polynomial $p(x)={\sum\limits_j} c_j x^j$ in
$n+1$ variables, where the arithmetic operations are
taken from the semiring $\R_t$. Depending on $t$ this
polynomial defines different functions $p_t:\R^{n+1}\to\R$.
Note that the function $$f_t(z)=t^{p_t(\Log_t(z))}$$
coincides with the patchworking polynomials where all $a_j=1$
and $v(j)=c_j$.
Here $\Log_t(z_1,\dots,z_{n+1})=(\log_t(z_1),\dots,\log_t(z_{n+1})).$

\begin{lem}
\label{normy}
If a point $x\in\R^{n+1}$ belongs to the amoeba
$$\Log_t(\{z\in\tor\ |\ f_t(z)=0\})$$ then
the monomials $c_jx^j$ from $p_t$
satisfy the generalized triangle inequality in $\R_t$, i.e.
for each index $k$ we have
$$c_k\odot x^k\le {\bigoplus\limits_{j\neq k}} c_j\odot x^j.$$
\end{lem}
\begin{proof}
If $x=\Log_t(z)$ with $f_t(z)=0$ then the sum of the monomials
$t^{c_j}z^j$ is zero and thus their norms must satisfy the triangle
inequality.
\end{proof}

Let $f_t=\sum\limits_{j\in\Delta\cap\Z^{n+1}}a_jt^{-v(j)}z^j$ now be
a general patchworking polynomial. Denote $\V_t=\{f_t=0\}\subset\tor$.
The family $f_t$ can be treated as a single polynomial in $(K^*)^{n+1}$
(see \ref{naa}). It defines a hypersurface $\V_K\subset(K^*)^{n+1}$.
Recall that {\em the Hausdorff distance} between two closed subsets
$A,B\subset\R^{n+1}$ is the number $$\max\{\sup\limits_{a\in A} d(a,B),
\sup\limits_{b\in B} d(b,A)\},$$ where $d(a,B)$ is the Euclidean distance
between a point $a$ and a set $B$ in $\R^{n+1}$.
Denote $\am_t=\Log_t(\V_t)$
and $\am_K=\Log_K(\V_K)$.

\begin{cor}
\label{ak}
The amoebas $\am_t$ converge in the Hausdorff metric
to the non-Archimedian amoeba $\am_K$ when $t\to\infty$.
\end{cor}
\begin{proof}
Lemma \ref{normy} and the inequality \eqref{Md} imply that
$\am_t$ converge to a subset of $\am_K$.
Indeed, for each $t$ we can rewrite $|a_jt^{v(j)}z^j|$ as
$|t^{c_j}z^j|$, $c_j=v(j)+\log_t|a_j|$. Such a monomial
induces a linear function $c_j+jx$ in $\R^{n+1}$.
The inequalities
\begin{equation}
\label{okr}
c_k+kx\le\max\limits_{j\neq k}(c_j+jx)+\log_t(N),
\end{equation}
where $N+1$ is the number of monomials in $f_t$, cut out
a uniformly bounded neighborhood of $\am_K$ which contains
$\am_K$.

The limit of $\am_t$ cannot be any smaller than $\am_K$
by the following topological reason.
A component of the complement of the set described by the
inequalities \eqref{okr} is given by the inequality
$c_k+kx>\max\limits_{j\neq k}(c_j+jx)+\log_t(N)$.
By \cite{FPT} this component is contained in the
component of $\R^{n+1}\setminus\am_t$ corresponding to
the index $k$. Thus, different components of the set described
by \eqref{okr} must be contained in different components of
$\R^{n+1}\setminus\am_t$.
\end{proof}

This corollary can be strengthened to describe the limits
of the varieties $\V_t\subset\tor$ under the corresponding
renormalization of the norms of their points.
The description is in terms of the lifts of non-Archimedian
amoebas, see \ref{lifts}.
Let $H_t:\tor\to\tor$ be the transformation defined by
$$H_t(z_1,\dots,z_{n+1})=
(t^{-|z_1|}\frac{z_1}{|z_1|},\dots,t^{-|z_{n+1}|}\frac{z_{n+1}}{|z_{n+1}|}).$$
We have $\Log_t=\Log\circ H_t$.

\renewcommand{\thethm}{5}
\begin{thm}
The sets $H_t(\V_t)$ converge in the Hausdorff metric to
$W(\V_K)$ when $t\to\infty$.
\end{thm}
The proof is the same as the proof of Corollary \ref{ak}.
The only difference we have to make is to incorporate the
arguments of the monomials to the inequalities \eqref{okr}.

\subsection{Construction of the fibration $\lambda_t:\V_t\to\Pi$}
\label{constlambda}
Let $\Pi$ be a maximal dual $\Delta$-complex and
$v:\Delta\cap\Z^{n+1}\to\R$ be the function such that
$\Pi=\Pi_v$ as in Proposition \ref{obratno}.
It gives us a patchworking polynomial
$f_t=\sum\limits_{j\in\Delta\cap\Z^{n+1}}t^{-v(j)}z^j$.
As before we denote with $\V_t\subset\tor$ the zero set
of this polynomial.

We construct $\lambda_t:\V_t\to\Pi$
for a sufficiently large $t$ by gluing the
fibrations $\lambda_H$ from \ref{hyperplanes}.

To do it we construct a singular foliation $\F_\Pi$ in a
neighborhood $\NN\supset\Pi$.
By Proposition \ref{prim-uj} $\Pi$ can be locally identified
with $\Sigma_n$ by elements of $ASL_{n+1}(\Z)$.
Recall that an element $M\in ASL_{n+1}(\Z)$ is a rotation defined by
a unimodular integer $(n+1)\times(n+1)$-matrix $(m_{j,k})$
followed by a translation by $m=(m_1,\dots,m_{n+1})$ in $\R^{n+1}$.
This transformation of $\R^{n+1}$ lifts to $\tor$ as
$$H_M:z_j\mapsto e^{m_j} z_1^{m_{j,1}}\dots z_{n+1}^{m_{j,n+1}}.$$

We patch the foliations $\F$ constructed in \ref{hyperplanes} for
the primitive $n$-complex $\Sigma_n$.
Let $v_j\in\Pi$ be a vertex. By Proposition \ref{prim-uj}
there exists a neighborhood $U_j\ni v_j$ in $\Pi$ and
$M_j\in ASL_{n+1}(\Z)$ such that $M_j(U_j)$ is a neighborhood
of $0$ in $\Sigma_n$.
Let $N_j$ be a small neighborhood
of the closure of $M_j(U_j)$.

Consider the pull-back under $M_j$ of the foliation $\F$ constructed
in \ref{hyperplanes} restricted to $N_j$.
Note that $M_j^{-1}(N_j)$
cover $\Pi$.
The pull-back foliations at the overlaps
$M_j^{-1}(N_j)\cap M_k^{-1}(N_k)$ do agree in general.
Nevertheless, they have the same type of singularities at the same points
and their non-singular leaves are transverse to $\Pi$.
A partition of unity gives a foliation $\F_\Pi$ in
a neighborhood $\NN$ of $\Pi$. Note that we can ensure that
$\NN$ contains an $\epsilon$-neighborhood of $\Pi$ for some
$\epsilon>0$.
Following \ref{hyperplanes} we denote $\pi_{\F_\Pi}:\NN\to\Pi$
the projection along the leaves of $\F_\Pi$.

By Corollary \ref{ak} for a sufficiently large $t>0$
we have $\Log_t(V_t)\subset\NN,$ and we define
$$\lambda_t=\pi_{\F_\Pi}\circ\Log_t:V_t\to\Pi.$$

\subsection{Proof of
Theorems \ref{comp} and \ref{locthm}}
\label{thm42}
Here we prove that $V_t$ is non-singular and that
$\lambda_t$ satisfies to all hypotheses of
Theorem \ref{tor} for a large $t>0$.

Note that if $\Log_t z=x$ then $||t^{-v(j)}z^j||=t^{jx-v(j)}$,
where $jx\in\R$ stands for the scalar product.
Let $F\subset\Pi$ be an open $(n+2-k)$-cell.
\begin{lem}
There exists $k$ monomials $t^{-v(j_1)}z^{j_1}, \dots,
t^{-v(j_k)}z^{j_k}$ that dominate $f_t$
in a neighborhood of $F$. Namely, any other monomial evaluated
at a point near $F$ has a smaller order by $t$.
Furthermore, the hypersurface
$$\sum\limits_{m=1}^k t^{-v(j_m)}z^{j_m}=0$$
is isomorphic to the hyperplane
$z_1+\dots+z_{k-1}+1=0$ under the multiplicative change
of coordinates by an element of $ASL_{n+1}(\Z)$.
\end{lem}
\begin{proof}
This follows from the maximality of $\Pi$.
By Proposition \ref{podrazd} $F$ is dual to a $k$-dimensional
polyhedron from a subdivision of $\Delta$. Since $\Pi$ is
maximal, this polyhedron is the standard $(k-1)$-simplex
up to action of $ASL_{n+1}(\Z)$.
\end{proof}

This lemma implies that $\V_t$ is non-singular for large $t>0$.
Indeed, it is covered by a finite number of open sets and
in each set it is a small perturbation of the image of a hyperplane.
Furthermore, its compactification $V_t\subset\C T_\Delta$
is smooth and transverse to the coordinate hyperplanes as
the same reasoning with the terms of smaller order
applies to the affine charts of $\C T_\Delta$.

Our next step is to isotop $V_t$ over $N_j$ as in \ref{lochyp}.
Recall that $N_j$
was defined in \ref{constlambda} as
a small neighborhood of $\bar{U_j}\subset\Pi$ in $\R^{n+1}$.
Denote $$Q^n_j=M_j^{-1}(H_t(Q^n))\cap\Log_t^{-1}(N_j).$$
By the last conclusion of Proposition \ref{locprop} these manifolds coincide
over $N_j\cap N_k$ for $t>>0$. We set
$$Q_\Pi=\bigcup\limits_j Q^n_j.$$
Note that for $t>>0$ $\V_t$ is isotopic to $Q_\Pi$ by the same
isotopy as in the proof of Proposition \ref{locprop} since all
other monomials of $f_t$ have smaller order in $t$.
This proves Theorem \ref{locthm}.
As in Proposition \ref{locprop} the closure
$\bar{Q_\Pi}\subset\C T_\Delta$ is a smooth manifold.
Similarly, $V_t$ is isotopic to $\bar{Q_\Pi}$ in $\C T_\Delta$.

In the proof of Theorem \ref{tor} we may assume that $\V=\V_t$ since
its closure $V_t\subset\C T_\Delta$ is smooth and transverse
to the coordinate hyperplanes. Similarly, in the proof of Theorems
\ref{diff}, \ref{diffg} and \ref{comp} we may assume that $V=V_t$.
We define $\lambda^\circ:\V\to\Pi$ as a composition of the isotopy
$V\approx Q_\Pi$,
the map $\Log_t:\tor\to\R^{n+1}$ and the projection
$\pi_{\F_\Pi}:\NN\to\Pi$.
Note that the isotopy $V\approx Q_\Pi$ is a symplectomorphism
by the Moser trick. (By the Moser trick, see e.g. \cite{CS} any smooth
deformation of a symplectic structure on a simply-connected
manifolds is isomorphic to the original symplectic structure.)
Note also that the Moser trick can be done equivariantly
with respect to the complex conjugation if $V$ is defined over $\R$.

To define $\lambda:V\to\bar{\Pi}$
we compactify the previous construction by using $\bar{Q_\Pi}$
and the reparametrized moment map to $\Delta$ as in \ref{compactif}.
\begin{figure}[h]
\centerline{\psfig{figure=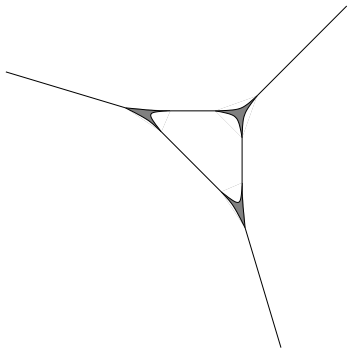,height=1.8in,width=1.8in}}
\caption{\label{locvt} The amoeba of the localization $Q_\Pi$
of a hypersurface.}
\end{figure}

This proves Theorem \ref{comp}, since everything in our construction
is equivariant with respect to complex conjugation as long as $a_j$
in the patchworking polynomial are real. The fibration $\lambda$
is totally real since it is totally real for a hyperplane.

Also, by Proposition \ref{locprop} this proves the second and the third
conclusions in Theorems \ref{diff}, \ref{diffg} and \ref{tor}.
The homotopy type of $\bar\Pi$ and $\Pi$ is the wedge
of $p_g$ copies of $S^n$, where $p_g=h^{n,0}$
by Proposition \ref{genus}.

To finish the proof of Theorems \ref{diff}, \ref{diffg} and \ref{tor}
we need to prove injectivity of the induced
homomorphism in cohomology and to exhibit the Lagrangian spheres
lifting the cycles from $\Pi$.

\subsection{Proof of Theorems \ref{diff}, \ref{diffg} and \ref{tor}}
The Lagrangian spheres will come from components
of certain real hypersurfaces whose complexification
is isotopic to $V$.


Let $j$ be a lattice point of $\Delta$.
We define $$f_t^{(j)}=\sum\limits_{k\neq j} |a_k|t^{v(k)}z^k-
|a_j|t^{v(j)}z^j.$$
Denote with $V_t^{(j)}\subset\tor$ the zero set of $f_t^{(j)}$
and with $\R V_t^{(j)}\subset\rtor$ its real part.
The Viro patchworking theorem
\cite{Vi} (see also \cite{GKZ} for a special case of {\em combinatorial
patchworking} and \cite{IV} for an elementary description in the case
of curves)
implies that $\R V_t^{(j)}\cap \R_+^{n+1}$
is diffeomorphic to a sphere $S^n$.
This sphere $S^n_j\subset V_t^{(j)}$
is Lagrangian
as a component of the real part
and it maps under $\Log_t$ to $\NN\supset\Pi$ for $t>>0$.
Furthermore, it realizes in $H_n(\Pi)$ the class corresponding to $j$
according to Proposition \ref{genus}.
\begin{figure}[h]
\centerline{\psfig{figure=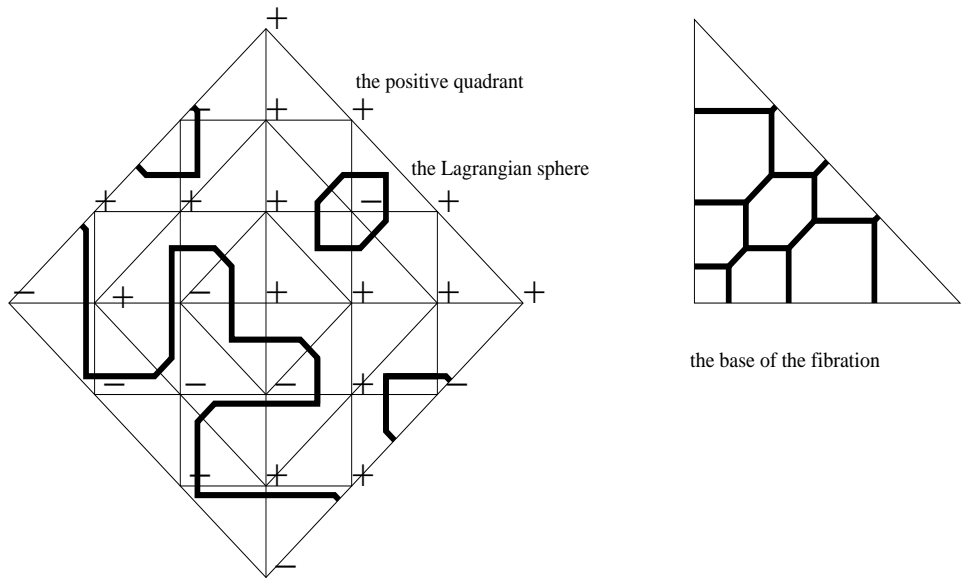,height=3in,width=5in}}
\caption{\label{patch} Construction of the Lagrangian
lift of a base cycle by the real patchworking}
\end{figure}

By \ref{thm42} $V_t^{(j)}$ is smooth. Thus,
it is isotopic to $V_t$ and we have
a diffeomorphism $h:V_t^{(j)}\to V_t$.
Moreover, we can choose an isotopy
among the hypersurfaces defined by the polynomials such that
the norm of all monomials is constant in the course of
deformation. All such hypersurfaces are smooth and their
image under $\Log_t$ is contained in $\NN\supset\Pi$ by
\ref{thm42}. Therefore, the image $h(S^n_j)$
projects to the same class in $H_n(\Pi)$.

By Moser's trick, $h$ is isotopic to
a symplectomorphism. This gives a Lagrangian sphere in
$V_t$ which projects to the class in $H_n(\Pi)$ corresponding
to $j$. Thus the last conclusion
of Theorems \ref{diff} and \ref{diffg} is proved.

Existence of such spheres also implies the first conclusion
of Theorems \ref{diff} and \ref{diffg}. The map $\lambda^*$
is injective since we can distinguish the images in $H^n(V;\Z)$
by their evaluations on these Lagrangian spheres.

The proof of Theorem 3 is the same since these
spheres belong to the toric part $\R\V_t$ of $\R V$.


\end{document}